\def\Date{July 7, 2010}

\magnification=\magstephalf
\hsize=14true cm
\vsize=24true cm
\frenchspacing\parskip4pt plus 1pt

\newread\aux\immediate\openin\aux=\jobname.aux 
\ifeof\aux\message{ <<< Run TeX a second time >>> }
\batchmode\else\input\jobname.aux\fi\closein\aux

\input suka.def\immediate\openout\aux=\jobname.aux

\draftfalse  
\ifdraft\footline={\hss\sevenrm\Date\hss}\else\nopagenumbers\fi

\headline={\ifnum\pageno>1\eightrm\ifodd\pageno\hfill~~ and tube
realizations of hyperquadrics
\hfill{\tenbf\folio}\else{\tenbf\folio}\hfill Classification of
commutative algebras \hfill\fi\else\hss\fi}

\centerline{\Gross Classification of commutative algebras and}
\medskip \centerline{\Gross tube realizations of hyperquadrics}

\bigskip\bigskip \centerline{By {\sl Gregor Fels} and {\sl Wilhelm
Kaup}} {\parindent0pt\footnote{}{\ninerm 2000 Mathematics Subject
Classification: 32V30, 13C05.}}

{\bigskip\medskip\narrower\ninerm\noindent\advance\baselineskip-1pt
ABSTRACT: In this paper we classify up to affine equivalence all local
tube realizations of real hyperquadrics in $\CC^{n}$. We show that
this problem can be reduced to the classification, up to isomorphism,
of commutative nilpotent real and complex algebras.  We also develop
some structure theory for commutative nilpotent algebras over
arbitrary fields of characteristic zero.\par}

\KAP{Introduction}{Introduction}

It is a well-known fact that every real-analytic manifold $M$ together
with an involutive CR-structure $(HM,J)$ admits at least locally a
generic embedding into some $\CC^n$, such that the CR-structure
induced from the ambient space $\CC^n$ coincides with the original
one.  A particularly important class of CR-submanifolds of $\CC^n$ are
the so-called CR-tubes, i.e., product manifolds $M=iF+\RR^n\subset
i\RR^n\oplus \RR^n=\CC^n$ together with the inherited CR-structure,
where $F\subset \RR^n$ is a submanifold. One important point here is
that the CR-structure of $iF+\RR^n$ is closely related to
real-geometric properties of the base $F$, which are often easier to
deal with, see e.g. \Lit{FKAU}.  In general, a CR-manifold will not
admit a local realization in $\CC^n$ as a CR-tube. On the other hand,
as shown by the example of the sphere $S=\{z\in \CC^n:|z_1|^2+\cdots
+|z_n|^2=1\}$, it is not immediate, that $S$ does admit several
affinely inequivalent local tube realizations, see \Lit{DAYA}.

It is quite obvious that the existence of a CR-tube realization for a
given CR-manifold $(M,HM,J)$ is related to the presence of certain
abelian subalgebras $\7v$ in $\hol(M)$, the Lie algebra of
infinitesimal CR-transformations, which are induced by all real
translations $z\mapsto z+x$, $x\in \RR^n$. It is perhaps a little bit
more subtle to give sufficient and necessary conditions for abelian
Lie subalgebras of $\hol(M)$ to give local CR-tube realization of $M$.
This has been worked out in \Lit{FKAP}. For short, let us call every
such abelian subalgebra a `qualifying' subalgebra of $\hol(M)$.
Curiously, the notion of locally affine equivalence among various
(germs of) tube realizations for a given $M$ proved to be less
appropriate for the study of CR-manifolds as it is too fine for many
applications: Even a homogeneous CR-manifold may admit an affinely
non-homogeneous tube realization, and in such a case the
aforementioned equivalence relation will give rise to uncountable many
equivalence classes of tube realizations.  A coarser equivalence
relation has been introduced in \Lit{FKAP} which seems to be most
natural in the context of CR-tubes.  Moreover, it is quite surprising
that under certain assumptions the pure geometric question of globally
affine equivalence can be reduced to the purely algebraic problem of
classifying conjugacy classes of certain maximal abelian subalgebras
of $\hol(M)$ with respect to a well-chosen group $G$.

The purpose of this paper is to give a full classification of all local
CR-tube realizations of every hyperquadric
$$S_{p,q}=\big\{[z]\in \PP(\CC^{m}):
|z_1|^2+\cdots+|z_p|^2=|z_{p+1}|^2+\cdots +|z_{m}|^2\big\}$$ in the
complex projective space $\PP(\CC^{m})$, where $m:=p+q$ and $p,q\ge1$,
applying the general methods from \Lit{FKAP}. The (compact)
hyperquadric $S_{p,q}$ is the unique closed orbit of
$\SU(p,q)\subset\SL(m,\CC)$ acting by biholomorphic transformations
on $\PP(\CC^{m})$. In this situation
$\hol(S_{p,q})=\hol(S_{p,q},a)\cong \su(p,q)$ holds for every $a\in
S_{p,q}$.  It is well known and easy to see that $S_{p,q}$ is locally
CR-equivalent to the affine real quadric in $\CC^{r}$, $r:=m-1$,
$$\Big\{z\in\CC^{r}:\Im(z_{r})=\sum_{1\le k<p}|z_{k}|^{2}\;\;
-\sum_{p\le k<r}|z_{k}|^{2} \Big\}$$ with non-degenerate Levi form of
type $(p{-}1,q{-}1)$. Therefore the classification
problem for local tube realizations for both classes is the same,
compare also \Lit{DAYA}, \Lit{FKAP}, \Lit{ISMI}, \Lit{ISAE},
\Lit{ISAV}, \Lit{ISEV}, \Lit{TANA} for partial results in this
context.

We have shown in \Lit{FKAP} that every abelian subalgebra $\7v\subset
\hol(M)$, which yields a tube realization, determines an involution
$\tau_{\7v}:\hol(M)\to \hol(M)$. For a given hyperquadric $S_{p,q}$
however, it turns out that all arising involutions $\tau_{\7v}$ are
conjugate in $\7g:=\hol(S_{p,q})\cong \su(p,q)$. We therefore fix an
involution $\tau:\7g\to \7g$ (with fixed point set $\7g^\tau \cong
\so(p,q)$) once and for all and reduce the classification of tube
realizations to the algebraic classification of all maximal abelian
subalgebras $\7v$ of $\7g$ contained in the $(-1)$-eigenspace of the
non-riemannian symmetric pair $(\7g,\tau)$ up to conjugation by
$\SU(p,q)$ (in fact, up to conjugation by the normalizer $G$ of
$\SU(p,q)$ in $\SL(m,\CC)$, but these two groups differ only if $p=q$,
and in this case the classification with respect to one group can
easily be derived from the classification up to conjugation with
respect to the other).  In contrary to the special case of toral
maximal subalgebras (i.e., Cartan subalgebras) $\7t\subset \7g$ only
little is known about the general case of arbitrary abelian maximal
subalgebras $\7v\subset \7g$.  The key point here is that after some
reduction procedures the conjugacy class of a maximal abelian
subalgebra $\7v\subset \7g^{-\tau}$ is completely determined by its
$D$-invariant (which is a first rough invariant of $\7v$, determined
by its toral part, see \Ruf{DF} for more details) and a finite set of
maximal abelian subalgebras $\7n_{\!\jmath}$, consisting of
ad-nilpotent elements only in $\su(p_\jmath,q_\jmath)$ and
$\7{sl}(m_{\!\jmath},\CC)$. Hence, the classification task reduces
essentially to the classification of ad-nilpotent abelian subalgebras
$\7n_j$ up to conjugation in $\SU(p_\jmath, q_\jmath)$,
resp. $\SL(m_\jmath,\CC)$. By our constructions, to every such
$\7n_{\!\jmath}$ there is associated a finite-dimensional commutative
associative nilpotent algebra $\5N_\jmath$ over $\FF=\RR$,
resp. $\FF=\CC$. Our main algebraic result is then the following

\Theorem{} Let $G=\SU(p,q)$ or $G=\SL(m,\CC)$ and $\tau:G\to
G$ an involutive automorphism with $G^\tau\cong \SO(p,q)$,
resp. $G^\tau\cong \SO(m,\CC)$. For any two maximal abelian
ad-nilpotent subalgebras $\7n_1, \7n_2\subset \7g$, contained in the
$(-1)$-eigenspace of $\tau$, the following conditions are equivalent:
\0$\7n_1$ and $\7n_2$ are conjugate by an element in $G$,
\1$\7n_1$ and $\7n_2$ are conjugate by an element in $G^\tau$,
\1the associated algebras $\5N_1$ and $\5N_2$ are isomorphic
as abstract $\FF$-algebras.\Formend

\medskip\noindent The commutative nilpotent algebras $\5N$ occurring
in the above theorem all have a 1-di\-men\-sional annihilator.  On the
other hand, given any nilpotent commutative $\FF$-algebra $\5N$ with
1-dimensional annihilator $\5A$, we construct an invariant of $\5N$
which is a certain non-degenerate symmetric 2-form
$b_\pi:\5N/\5A\times \5N/\5A\to \5A$. Depending on the type of $b_\pi$,
the algebra $\5N$ gives rise to a maximal abelian subalgebra in
$\7g^{-\tau}$ ($\7g$ and $\tau$ as in the preceeding theorem) and in
turn to a tube realization of a hyperquadric.

Summarizing, the classification of all local tube realizations of
hyperquadrics is essentially equivalent to the classification of
finite-dimensional nilpotent commutative $\FF$-algebras up to
isomorphism.  In this paper we give an explicit classification for low
values of $p$ or $q$, i.e., we carry out all those cases $(p,q)$ where
there are only finitely many isomorphy classes.  The classification in
terms of explicit lists seems to be hopeless in the general case. For
big values of $p$ and $q$ there are always uncountable many
inequivalent nilpotent commutative $\FF$-algebras.

As the algebraic results developed in this paper might be of broader
interest, we collect in the Appendix all relevant results
concerning the fine structure of nilpotent commutative algebras. These
are formulated in a more general setup (e.g. over arbitrary fields of
characteristic zero).

For certain applications one would like to have explicit defining
equations for the various tube realizations, determined by qualifying
subalgebras $\7v\subset \su(p,q)$. One of our main geometric results
is a procedure which produces for every qualifying $\7v$ an explicit
defining equation which describes the tube realization $iF_{\7v}\oplus
V_{\7v}$ of $S_{p,q}$.  In that way we obtain quite transparent
formulae, reflecting the algebraic structure of $\7v$.

\medskip The paper is organized as follows: In Section 2 we relate our
results to existing results in the literature, in particular to those
in \Lit{DAYA}, \Lit{ISMI}, \Lit{ISAE}. In Section 3 we recall the
necessary tools from \Lit{FKAP} and give a short outline of the
classification procedure. In particular we introduce certain abelian
subalgebras $\7v\subset\su(p,q)$ as {\sl qualifying MASAs} -- these
are the algebraic objects to be classified. In Section 4 we split
every MASA $\7v$ into its toral $\7v\red$ and its nilpotent part
$\7v\nil$ and classify the centralizers of $\7v\red$. Crucial for the
classification is the decomposition given by Lemma 4.6 that leads to a
combinatorial invariant $\3D(\7v)$ that we call the $D$-invariant of
$\7v$. For fixed $p,q$ the set $\5D_{p,q}$ of all $D$-invariants of
MASAs in $\su(p,q)$ is finite, but still, in general there are
infinitely many equivalence classes of MASAs in $\su(p,q)$ with a
fixed $D$-invariant. In Sections 5 and 6 we study MANSAs (maximal
commutative associative nilpotent subalgebras) in $\su(p_{j},q_{j})$
and $\7{sl(m_{j},\CC)}$ as these are the building blocks for general
MASAs in $\su(p,q)$. In Section 7 we demonstrate briefly how for every
MANSA $\7v\subset\su(p,q)$ with corresponding tube realization
$iF+\RR^{n}\subset\CC^{n}$ of $S_{p,q}$ the base $F\subset\RR^{n}$ can
be written in terms of a canonical equation. In section 9 we give two
examples of MANSAs and in the Appendix 10 we collect several algebraic
tools needed in the paper that might also be of independent interest.

\KAP{Preliminaries}{Preliminaries} In the following we characterize
algebraically the local tube realizations of the hyperquadric
$S=S_{p,q}\subset\PP_{r}:=\PP(\CC^{r+1})$ with $p,q\ge1$ and
$r:=p+q-1$, compare (3.3) in \Lit{FKAP}.  Since $S_{p,q}$ and
$S_{q,p}$ only differ by a biholomorphic automorphism of the
projective space $\PP_{r}$ it would be enough to discuss the case
$p\ge q$.

 The local tube realizations of $S$ up to affine equivalence in the
cases $q=1,2$ were obtained in the papers \Lit{DAYA}, \Lit{ISMI}
respectively by solving certain systems of partial differential
equations coming from the Chern-Moser theory \Lit{CHMO}. A
classification of the case $q=3$ has been announced in \Lit{ISAE},
proofs are intended to appear in the forthcoming book \Lit{ISEV}.

In this paper we give a classification for arbitrary $p,q$. It turns
out that this, after several reducing steps, essentially boils down to
the classification of abstract abelian nilpotent real and complex
algebras $\5N$ of dimension $r:=p+q-1$ with 1-dimensional annihilator.
For small values of $q$ these can be determined explicitly while for
large $p,q$ this appears to be hopeless.  On the other hand, we
associate to every local tube realization of $S_{p,q}$ a combinatorial
invariant $\3D$ out of a finite set $\5D_{p,q}$ in such a way that for
any two local tube realizations the equations for the corresponding
tube bases $F,\tilde F\subset\RR^{r}$ are essentially of the same type
(up to some polynomial terms in the coordinates of $\RR^{r}$ coming
from the aforementioned different abelian nilpotent algebras arising
naturally in this context).

To compare this with the known results in case $q\le3$ let us
introduce the number $n:=r-1=p+q-2$, so that every local tube
realization $T$ of $S=S_{p,q}$ is a hypersurface in $\CC^{n+1}$ with
CR-dimension $n$, and $r$ is the rank of the Lie algebra $\hol(S)$.
Let furthermore $c_{p,q}$ be the cardinality of all affine equivalence
classes of closed tube submanifolds in $\CC^{r}$ that are locally
CR-equivalent to $S$.\nline In case $q=1$, that is the case of the
standard sphere in $\CC^{p}$, \Lit{DAYA} implies
$c_{p,1}=p+2=n+3$.\nline In case $p\ge q=2$ we have $n=p$ and the
explicit list of tube realizations in \Lit{ISMI} implies the estimate
$c_{p,2}\le p(p+9)/2$. Our considerations will give
$$c_{p,2}=5p+k(p-k)-\delta_{p,2}\steil{with}k:=\lceil
p/2\rceil\,,\Leqno{HB}$$ where for every $t\in\RR$ the {\sl ceiling}
$\lceil t\rceil$ is the smallest integer $\ge t$ and $\delta$ is the
Kronecker delta.  Therefore, the list in the Theorem of \Lit{ISMI}
\p442 must contain repetitions. Indeed, in type \kl7 for every $s$ the
parameters $t$ and $\tilde t:=n-2+s-t$ give affinely equivalent tube
realizations. The same holds in case $p=2$ for $s=1$ in type \kl1 and
$s=0$ in type \kl2.\nline In case $p\ge q=3$ it has been announced in
\Lit{ISAE} that $c_{p,3}$ is finite if and only if $p\le5$.

In case $p,q\ge4$ we show that $c_{p,q}$ always is infinite. Except
for $c_{4,4,}$ this has already been announced in \Lit{ISAE}.

\KAP{setup}{The algebraic setup}

For fixed $p,q\ge1$ with $m:=p+q\ge3$ let $\2E\cong\CC^{m}$ be a
complex vector space and $\6h:\2E\times \2E\to\CC$ a hermitian form of
type $(p,q)$ ($p$ positive and $q$ negative eigenvalues). Since any
two hermitian forms of the same type on $\2E$ are equivalent (up to a
positive multiplicative constant) with respect to the group
$L:=\SL(\2E)\cong\SL(m,\CC)$ it does not matter which $\6h$ has been
chosen above. More important for computational purposes is to choose a
convenient vector basis of $\2E$ in such a way that the corresponding
matrix representation of $\6h$ is optimally adapted.

The complex Lie group $L$ acts in a canonical way transitively
on the complex projective space $Z:=\PP(\2E)$ with finite
kernel of ineffectivity (the center of $L$). The subgroup $$G:=\{g\in
L:\6h(gz,gz)=\pm\6h(z,z)\steil{for all}z\in\2E\}\Leqno{GP}$$ is a real
Lie group with $(1+\delta_{p,q})$ connected components acting
transitively on the hypersurface
$$S=S_{p,q}:=\{[z]\in\PP(\2E):\6h(z,z)=0\}\;.$$ 
We write for
the corresponding Lie algebras
$$\7l=\7{sl(\2E)}\Steil{and}\7g:=\su(\2E,\6h)=\{\xi\in\7l:\Re\6h(\xi
z,z)=0\steil{for all}z\in\2E\}\;.$$ As a matter of fact, $\7l$
coincides with the complex Lie algebra $\hol(Z)$ of holomorphic vector
fields on $Z$. Further, for every $a\in S$ the canonical inclusions
$\7g\into \hol(S)\into \hol(S,a)$ turn out to be isomorphisms, and
therefore we identify $\hol(S)$ with $\7g$.  With $\sigma:\7l\to\7l$
we denote the antilinear involutive Lie automorphism with
$\Fix(\sigma)=\7g$.

The hyperquadric $S=S_{p,q}$ satisfies the assumptions of Theorem 7.1
in \Lit{FKAP}, and the various tube realizations are, up to the global
affine equivalence as defined in \Lit{FKAP} Definition 6.1, in a
1-1-correspondence to $\Glob(S,a)$-conjugacy classes of certain
abelian subalgebras in $\hol(S,a)$ (see \Lit{FKAP} for the definition
of $\Glob(S,a)$ and its basic properties).  In the case under
consideration $\Glob(S,a)=\Ad(G)\subset\Aut(\hol(S,a))$ for every
$a\in S$; our task then will be to classify up to the action of
$\Ad(G)=\Ad(N_L(\7g))$ on $\7l$ all $\sigma$-invariant abelian
subalgebras $\7e\subset\7l$ which have an open orbit in $Z$.  Every
such $\7e$ automatically has complex dimension $r:=m{-}1$ and is
maximal abelian in $\7l$ by Lemma 2.1 in \Lit{FKAP}.

Every involution $\tau_{\7v}:(S,a)\to (S,a)$ extends to a global
involution $\widehat \tau_{\7v}:Z\to Z$. Moreover, any two such
involutions are conjugate by an element of $G$ (even of the connected
identity component $\SU(\2E,\6h)$ of $G$), compare \Lit{FKAP}.  The
search can therefore be restricted by fixing once and for all an
involution $\tau$ of $S$ whose fixed point set $S^{\tau}=\Fix(\tau)$
is not empty and has dimension $r-1$. Such a $\tau$ has a unique
extension to an antiholomorphic involution of $\PP(\2E)$ that comes
from a conjugation $\2E\to\2E$, $z\mapsto\overline z$, that is,
$\tau[z]=[\overline z]$ for all $[z]\in S$. By our results it is
enough to classify up to conjugation by $G$ all abelian Lie
subalgebras $\7v\subset\7g^{-\tau}$ with
$\epsilon_{a}(\7v)=\T^{-\tau}_{a}S$ for a given point $a\in S$. These
$\7v$ are automatically maximal abelian in $\7g\cong\su(p,q)$ and have
dimension $r=m{-}1=\rank(\7g)$.

\medskip\Joker{SE}{Setup} For the rest of the paper we fix the
following notation: For $p,q$ and $m=p+q$ as above, $\2E$ is a complex
vector space of dimension $m$ with (positive definite) inner product
$(z|w)$ (complex linear in the first and antilinear in the second
variable).  Furthermore, ${\tau:\2E\to\2E}$, $z\mapsto\overline z$, is
an (antilinear) conjugation on $\2E$ with $(\overline z|\overline
w)=(w|z)$ for all $z,w\in\2E$. With the same symbol $\tau$ we also
denote the induced antiholomorphic involution of $Z=\PP(\2E)$ as well
as of the complex Lie algebra $\7l:=\7{sl}(\2E)=\hol(\PP(\2E))$.  In
addition, $(e_{j})_{1\le j\le m}$ is an orthonormal basis of $\2E$ with
$e_{j}=\overline e_{j}$ for all $j$, and the hermitian form
$\6h=\6h_{p,q}$ on $\2E$ is given by
$$\6h(e_{j},e_{k})=\vartheta_{p,j}\delta_{j,k}
\steil{with}\vartheta_{p,j}:= \cases{\phantom{-}1&if $p\ge
j$\cr-1&otherwise .\cr}\Leqno{SB}$$ The involution $\tau$ on $Z$
leaves the hyperquadric $S=S_{p,q}$ invariant. Therefore also
$\7g=\su(p,q)=\hol(S)$ is invariant under the involution $\tau$ of
$\7l$.  As before, $\sigma$ is the involution of $\7l$ defining the
real form $\7g$ of $\7l$. Clearly, the involutions $\sigma$, $\tau$
commute on $\7l$.

$\SU(\2E,\6h)=\SU(p,q)$ is the
connected identity component $G^{0}$ of the group $G$ defined in \Ruf{GP}.
Only in case $p=q$ the group $G$ is disconnected and then
$\kmatrix{\;0}{\;\;\one}{-\one}0\in\SL(2p,\CC)$ is contained in the second
connected component of $G$.

With $\End(\2E)$ we denote the endomorphism algebra of $\2E$, a unital
complex associative algebra with involution $g\mapsto g^{*}$ (the
adjoint with respect to the inner product). With respect to the Lie
bracket $[f,g]=fg-gh$ it becomes a reductive complex Lie algebra that
is denoted by $\gl(\2E)$ and contains $\7{sl}(\2E)$ as semisimple
part. For every $z,w\in \2E$ we denote by $z\otimes w^{*}\in\gl(\2E)$
the endomorphism $x\mapsto(x|w)z$. Then $(z\otimes w^{*})^{*}=w\otimes
z^{*}$ is obvious. We also consider adjoints with respect to $\6h$ and
write $g^{\star}$ for the endomorphism satisfying
$\6h(gw,z)=\6h(w,g^{\star}\!z)$ for all $w,z\in \E$.

\bigskip

\Joker{TA}{The task} Let $G,\7g,\7l=\7{sl}(\2E),\sigma$ be as before
and let a compatible conjugation $\tau:\7l\to \7l$ be fixed once and
for all, induced by $z\mapsto \overline z$ on $\2E$.  In order to
classify all local tube realizations of $S=S_{p,q}$ up to globally
affine equivalence (compare Section 6 in \Lit{FKAP}) we have to
classify all abelian subalgebras
$\7v\subset\7g=\hol(S_{p,q})=\su(p,q)$ up to conjugation with respect
to $\Ad(G)$ which have the following property:

\item{(A)} The complexification $\7v^{\CC}$ has an open orbit in
$Z=\PP(\2E)$, that is, $\epsilon_{a}(\7v^{\CC})=\T_{a}Z$ for
some $a\in Z$ (and hence even for some $a\in S$).

\noindent This condition, justified by Proposition 4.2 in \Lit{FKAP}, 
is of geometric nature but implies the
following purely algebraic properties:

\item{(B)} $\7v$ is maximal abelian in $\7g$ -- we call every such
subalgebra a MASA in $\7g$.

\item{(C)} $\dim\7v=\rank\7g=\dim Z\quad (=r:=m-1)$.

\item{(D)} $\Ad(g)(\7v)\subset\7g^{-\tau}$ for some $g\in G$.

\noindent Instead of classifying all $G$-conjugation classes of $\7v$
with property (A) we classify more generally the classes of $\7v$
satisfying (B) and (D), let us call them {\sl qualifying MASAs} in
$\7g$ for the following. It will turn out a posteriori that these
$\7v$ automatically satisfy (A) and hence also (C).

\medskip\Joker{OL}{A short outline of the classification procedure} We
proceed by analyzing the algebraic structure of maximal abelian
subalgebras $\7v\subset\7g$.  We will need some well known facts from
the structure theory of semisimple Lie algebras (we refer to
\Lit{KNAP} and \Lit{WARN} as general references).  Write
$$\eqalign{N_{\7a}(\7b):&=\{x\in \7a:[x,\7b]\subset\7b\}\hbox{ for the
normalizer and}\cr C_{\7a}(\7b):&=\{x\in \7a:[x,\7b]=0\}\hbox{ \rm for
the centralizer}\cr}$$ of any subalgebra $\7b$ in a Lie algebra
$\7a$. Also let $Z(\7a):=C_{\7a}(\7a)$ be the center of $\7a$.  The
classification idea is based on the observation that each maximal
abelian $\7v\subset \7g\subset \End(\E)$ has a unique decomposition
into toral and nilpotent part, i.e., $\7v=\7v\red\oplus\7v\nil$, where
$\7v\red$ consists of semisimple and $\7v\nil$ of nilpotent elements
in $\End(\E)$. Each toral subalgebra, in particular $\7v\red$ of a
qualifying MASA $\7v$, gives rise to the real reductive subalgebra
$C_{\7g}(\7v\red)$.  On the other hand, the maximality of $\7v$
implies that $\7v\red=Z(C_{\7g}(\7v\red))$.  Hence, there is a natural
bijection between [the $G$-conjugacy classes of] toral parts $\7v\red$
of qualifying MASAs $\7v$ and [the $G$-conjugacy classes of] certain
reductive subalgebras $C_{\7g}(\7v\red)$. A particular class of
qualifying MASAs is formed by those $\7v$-s for which $\7v\nil=0$,
i.e., $\7v$ is a real Cartan subalgebra of $\7g$. It turns out that
each of the $\min\{p,q\}{+}1$ conjugacy classes (with respect to $G^0$
- or equivalently to $G$) of real CSAs has qualifying representatives.

It is well-known that general centralizers of tori in {\sl complex}
semisimple Lie algebras can be characterized by subsets of simple
roots.  In our case however we have to classify {\sl real}
centralizers. An additional complication is that not all (conjugacy
classes of) centralizers of tori are of the form $C_{\7g}(\7v\red)$
with a qualifying MASA $\7v$.  In the first part of our classification
we provide a combinatorial tool giving an explicit characterization of
these conjugacy classes of centralizes which are related to qualifying
MASAs.

The nilpotent part $\7v\nil$ of a qualifying MASA is contained in the
semisimple part of $C_{\7g}(\7v\red)$, more precisely, we have the
following diagram:
$$\diagram{ \7v & \;=\; & \7v\red & \oplus & \7v\nil \cr \cap & &
\Vert \quad & & \cap\quad \cr C_{\7g}(\7v\red)^{-\tau}&=&
Z(C_{\7g}(\7v\red))& \oplus & C_{\7g}\ss(\7v\red)^{-\tau} \rlap{
.}}\qquad \Leqno{MD}$$ Moreover, $\7v\nil$ is a maximal abelian and
nilpotent subalgebra of $C_{\7g}\ss(\7v\red)$ -- we call such
subalgebras MANSAs.  An important observation is that the simple
factors occurring in $C_{\7g}\ss(\7v\red)$ are not arbitrary real
forms of $C_{\7l}(\7v\red)$: A simple factor in $C_{\7g}\ss(\7v\red)$
is isomorphic either to $\su(p',q')$ or ${\7{sl}}_m(\CC)$.
Consequently, $\7v\nil$ is a product of qualifying MANSAs in
$\su(p',q')$ and ${\7{sl}}_m(\CC)$.  The classification of the last
mentioned Lie subalgebras turns out to be equivalent to the
classification of arbitrary real or complex {\sl associative}
commutative nilpotent algebras $\5N$ with 1-dimensional annihilator.

 We will analyze the consequences of condition (C) later on; one
outcome is that $\dim \7v\nil=\rank(C_{\7g}\ss(\7v\red))$, hence $\dim
\7v=\rank \7g$ and each such $\7v^\CC$ has an open orbit in $Z$.
Summarizing, our task then is reduced to the solution of the following
algebraic problems: \item{[R]} Classify up to conjugation all
reductive subalgebras $\7r\subset \7g$ which are centralizers of the
reductive part of a qualifying MASA $\7v\subset \7g$ (compare
\ruf{TA}).  \item{[N]} Given a $\tau$-stable reductive subalgebra
$\7r=Z(\7r)\oplus \7r\ss$ of the above type, classify up to
conjugation all maximal abelian nilpotent subalgebras $\7n$ of
$\7r\ss$ with $\7n \subset (\7r\ss)^{-\tau}$.

\KAP{centralizers}{Classification of the centralizers 
$C_{\7g}\big(\7v\red\big)$}

\Joker{HA}{Qualifying Cartan subalgebras.} Particular examples of MASAs
$\7v\subset \7g$ are maximal toral subalgebras, i.e., Cartan
subalgebras.  This is precisely the case when $\7v=\7v\red$ and
$\7v\nil=0$. It is well-known that $\su(p,q)$ (with $q\le p$) has
$q+1\;$ $G^0$--conjugacy classes of real Cartan subalgebras.  We need
to know that each such conjugacy class contains a {\sl qualifying}
MASA of $\7g$:

\medskip \Lemma{CT}{\bf (Maximal toral subalgebras)} Let $\7g$, $\tau$
and $S=S_{p,q}$ with $p\ge q$ be as before. Then: \0 Every complex
Cartan subalgebra of $\7{sl}(\2E)$ has an open orbit in $\PP(\2E)$.  \1
Every $G^{0}$-conjugacy class of real Cartan subalgebras in $\7g$ has a
representative contained in $\7g^{-\tau}$.  \1 To every maximal
abelian subalgebra $\7v\subset \7g^{-\tau}$ there exists a real Cartan
subalgebra $\7h$ of $\7g$ and a $g\in G$ with
$\Ad(g)(\7v\red)\subset\7h\subset\7g^{-\tau}$.

\Proof (i) is an easy consequence of the fact that there is only one
conjugacy class of complex CSAs in $\7{sl}(\2E)$ and that the subspace
of all diagonal matrices in $\7{sl}(\2E)$ is one of them. For the proof
of (ii) fix an arbitrary real CSA $\7h$ of $\7g$. Then $\7h^{\CC}$ is
a complex CSA of $\7{sl}(E)$ and hence has an open orbit in
$\PP(\2E)$. Therefore, by Propositions 4.2 and 3.2 in
\Lit{FKAP}, there is a point $a\in S$ and an involution $\theta$ of
$(S,a)$ with $\7h\subset\7g^{-\theta}$. Also, $\theta$ satisfies
(3.1) in \Lit{FKAP} and extends to an antiholomorphic involution of
$\PP(\2E)$. Therefore $\tau=g\theta g^{-1}$ for some $g\in G^{0}$, that
is, $\Ad(g)(\7h)\subset\7g^{-\tau}$. Below, we also give an
alternative, algebraic proof of (ii) without referring to results from
\Lit{FKAP}.  Assertion (iii) follows from (ii) since there exists a
CSA $\7h$ of $\7g$ with $\7v\red\subset\7h$.\qed

Every toral subalgebra $\7t\subset\7g\subset\End(\2E)$ has a unique
decomposition $\7t=\7t_{+}\!\oplus\,\7t_{-}$ into its compact and its
vector part, that is, all elements in $\7t_{+}$ ($\7t_{-}$
respectively) have imaginary (real respectively) spectrum as operators
on $\2E$. Clearly the dimensions of these parts are invariants of the
$G$-conjugacy class of $\7t$ in $\7g$, and in the case of a semisimple
Lie algebra $\7g$ of Hermitian type, (as for instance $\su(p,q)$)
$\dim \7t_+$ determines uniquely its conjugacy class.  For later use
we construct explicitly for every $\ell=0,1,\dots,q$ a CSA
$^{\ell}\!\7h$ of $\7g$ with $^{\ell}\!\7h\subset\7g^{-\tau}$ and
$\ell=\dim(^{\ell}\7h_{-})\,$.

\Joker{}{Diagonal bases.} Consider on the integer interval
$\{1,2,\dots,m\}$ the reflection defined by $j\mapsto
j^{\bullet}:=m{+}1{-}j$ and recall the choice of the 
orthonormal basis $(e_{j})_{1\le j\le m}$ and of $\vartheta_{p,j}$
as in \Ruf{SB}.  Fix an integer $\ell$ with $0\le\ell\le q$, a
complex number $\omega$ with $2\omega^{2}=i$ and define a new
orthonormal basis $(^{\ell}\!f_{j})_{1\le j\le m}$ of $\2E$ by
$$^{\ell}\!f_{j}:=\cases{ e_{j}&if $\ell<j<\ell^{\bullet}$\cr \omega
e_{j}+\overline\omega e_{j^{\bullet}}& otherwise\cr },\steil{for
which} e_{j}=\cases{ ^{\ell}\!f_{j}&if
$\ell<j<\ell^{\bullet}$\vadjust{\vskip3pt}\cr \overline\omega\;
^{\ell}\!f_{j}+\omega\; ^{\ell}\!f_{j^{\bullet}}& otherwise\cr }$$ is
easily verified. Then for all $1\le j\le k\le m$ we have
$$ \6h(^{\ell}\!f_{j},^{\ell}\!f_{k})=\cases{ i\delta_{j,k^\bullet}&
if $j\le \ell$ \cr \vartheta_{p,j}\delta_{j,k}&if $\ell< j<
\ell^\bullet$ .\cr }\Leqno{LF}$$
Let $^{\ell}\!\7h\subset\7g=\su(\2E,\6h)$ be the abelian
subalgebra of all endomorphisms that are diagonal with respect to
$(^{\ell}\!f_{j})$. From

$$\eqalign{
{}^{\ell}\!f_{j}\otimes\,^{\ell}\!f_{j}^{*}\,+\,^{\ell}\!f_{j^{\bullet}}
\otimes\,^{\ell}\!f_{j^{\bullet}}^{*} &=e_{j}\otimes\,e_{j}^{*}\,+
\,e_{j^{\bullet}}\otimes\,e_{j^{\bullet}}^{*}\cr
{}^{\ell}\!f_{j}\otimes\,^{\ell}\!f_{j}^{*}\,-\,^{\ell}\!f_{j^{\bullet}}
\otimes\,^{\ell}\!f_{j^{\bullet}}^{*}
&=ie_{j}\otimes\,e_{j^{\bullet}}^{*}\,-
ie_{j^{\bullet}}\otimes\,e_{j}^{*}\;.\cr }$$ for all $j\le\ell$ we
derive that the decomposition
${^{\ell}\!\7h=^{\ell}\!\7h_{+}\oplus\;^{\ell}\!\7h_{-}}$ into compact
and vector parts is given by
$$\eqalign{ {^\ell}\!\7h_{-}&=\bigoplus_{j\le\ell} i\RR\Big(
e_{j}\otimes\,e_{j^{\bullet}}^{*}\,- e_{j^{\bullet}}\otimes\,e_{j}^{*}
\Big)\cr
{^\ell}\!\7h_{+}&=\Big(\bigoplus_{j\le\ell}i\RR\big(e_{j}\otimes\,e_{j}^{*}\,+
\,e_{j^{\bullet}}\otimes\,e_{j^{\bullet}}^{*}\big)
\;\oplus\!\!\bigoplus_{\ell<j<\ell^{\bullet}}\!\!\!
i\RR(e_{j}\otimes\,e_{j}^{*})\Big)_{\hbox{\rm tr}=0}\;.\cr }\Leqno{ZS}$$ As a
consequence, ${}^{\ell}\!\7h=C_{\7g}({}^{\ell}\!\7h)$  is a CSA of
$\7g$ with $\dim({}^{\ell}\7h_{-})=\ell$ and 
${}^{\ell}\!\7h\subset\7g^{-\tau}$. This gives a constructive
proof for (ii) in Lemma \ruf{CT}.\qed

\Joker{}{The general case $\7v\red\subset\7v$.}
We proceed to the general case where $C_{\7g}(\7v\red)$ may contain
$\7v\red$ properly, that is $\7v\nil\ne0$. Given $\7v\red$, or
equivalently $C_{\7g}(\7v\red)$, after conjugating with an element of
$\SU(\2E,\6h)$ we may assume that
$\7v\red\subset{}^{\ell}\7h\subset\7g^{-\tau}$ for some $\ell\le q$ as
above.
In the complex situation, (i.e., for the centralizer in $\7g^\CC$, or
equivalently in $\gl(\2E)$) it is well known that there is a unique
direct sum decomposition with summands $\2E_{\jmath}\ne0$
$$\2E=\bigoplus_{\jmath\in\5J}\2E_\jmath\Steil{such that}
C_{\gl(\2E)}(\7v\red)=\bigoplus_{\jmath\in\5J}
\7{gl}(\2E_\jmath)\,.\Leqno{ZR}$$ The subspaces $\E_j$ correspond to
joint eigenspaces of the toral abelian subalgebra $\7v\red\subset
\7{sl}(\E)$ with respect to certain functionals $\gamma\in
(\7v\red)^*$. Since $\7v\red$, and in turn $C_{\7{l}}(\7v\red)$ is
invariant under the conjugation $\tau$, we conclude that there is an
involution $\jmath\mapsto\overline \jmath$ of the index set $\5J$ with
$\tau\2E_{\jmath}=\2E_{\bar \jmath}$ for all $\jmath\in\5J$.  One key
point here is that for every $\jmath\in\5J$ the restriction
$\6h_{\jmath}$ of $\6h$ to $\2E_{\jmath}+\2E_{\bar \jmath}$ is
non-degenerate while in case $\bar\jmath\ne\jmath$ the spaces
$\2E_{\jmath}$ and $\2E_{\bar \jmath}$ are totally $\6h$--isotropic
and have zero intersection.  (A priori, the non-degeneracy of the
restrictions $\6h_\jmath$ does {\sl not} follow from the mere
$\tau$--invariance of the decomposition $\2E=\bigoplus
\2E_\jmath$. However, since $\7v\red\subset {}^\ell\!\7h$ for some
$\ell$, one can show using root theory that every subspace
$\2E_\jmath\subset\2E$ occurring in the above decomposition is
invariant under every orthogonal projection
${}^\ell\!f_k\otimes{}^\ell\!f_k^{*}$, $k\in\5J$. With \ruf{LF} then
the non-degeneracy of $\6h_{\jmath}$ follows).

\Joker{}{Restrictions of $\sigma$ and $\tau$ to the simple factors of
the centralizer.}  Choose a subset $\5L\subset\5J$ such that
$\5J=\5K\cup\5L\cup\overline{\5L}$ is a disjoint union for
$\5K:=\{\jmath\in\5J:\bar\jmath=\jmath\}$.  For every $\jmath\in\5K$
the subalgebra $\7{sl}(\2E_{\jmath})\subset\7{sl}(\2E)$ is invariant
under $\tau$ as well as $\sigma$, that is
$$\7{sl}(\2E_{\jmath})^{\sigma}= \su(\2E_{\jmath},h_{\jmath})\qquad
\7{sl}(\2E_{\jmath})^{\tau}=\7{sl}(\2E_{\jmath}^\tau)$$ (with
$\su(\2E_{\jmath},h_{\jmath})=0=\7{sl}(\2E_{\jmath})$ in case
$\dim\E_{\jmath}=1$). Also, for every $\jmath\in\5L$ we have
$\E_{\bar\jmath}=\tau(\E_\jmath)$, $
\sigma(\7{sl}(\E_{\jmath}))=\tau(\7{sl}(\E_{\jmath}))=\7{sl}(\E_{\bar\jmath})$
and
$\sigma,\tau\in\Aut_\RR(\7{sl}(\E_{\jmath})\oplus\7{sl}(\E_{\bar\jmath}))$
are given by
$$\tau(x,y)=(\underline\tau y\underline \tau,\underline\tau
x\underline\tau)\Steil{and}\sigma(x,y)=(-y^{\star},-x^{\star})\;,$$ where
$x^{\star},y^{\star}$ are the adjoints of $x,y$ with respect to the
hermitian form $\6h$. The symmetric complex bilinear form
$\beta:\E_{\jmath}\times\E_{\jmath}\to\CC$ defined by
$\beta(x,y)=\6h(x,\tau y)$ is non-degenerate and
$$\big(\7{sl}(\E_{\jmath})\oplus\7{sl}(\E_{\bar\jmath})\big)^{\sigma}\cong
\5R^{\CC}_{\RR}\big(\7{sl}(\E_{\jmath})\big)\quad\big(\7{sl}(\E_{\jmath})
\oplus\7{sl}(\E_{\bar\jmath})\big)^{\sigma\tau}\cong
\5R^{\CC}_{\RR}\big(\so(\E_{m_{\jmath}},\CC)\big)\,,\Leqno{NT}$$ where
$\5R_{\RR}^{\CC}$ is the forgetful functor restricting scalars from
$\CC$ to $\RR$.  With these ingredients we can state:

\Lemma{EI} For the decomposition \Ruf{ZR} we have
$$ \eqalign{ \7v\red=Z(C_{\7g}(\7v\red))&= \Big(
\bigoplus_{\5K\cup\5L}i\RR\id_{\2E_{\jmath}+\2E_{\bar\jmath}}\Big)_{\tr=0}\;\oplus\;
\bigoplus_{\5L}\RR(\id_{\2E_{\jmath}}-\id_{\2E_{\bar\jmath}}) \cr
C_{\7g}\ss(\7v\red)&=\;\bigoplus_{\5K} \su(\2E_{\jmath},\6h_{\jmath})
\oplus \bigoplus_{\5L} \big(\7{sl}(\2E_\jmath)\oplus
\7{sl}(\2E_{\bar\jmath}) \big)^\sigma\cr&\;\cong\bigoplus_{\5K}
\su(p_{\jmath},q_{\jmath}) \oplus\;
\bigoplus_{\5L}\7{sl}(m_{\jmath},\CC)\;,}$$ where
$m_{\jmath}=\dim(\2E_{\jmath})$ and $(p_{\jmath},q_{\jmath})$ for
$\jmath\in\5K$ is the type of the restriction $\6h_{\jmath}$ on
$\2E_{\jmath}$.  If $\7v\red\subset\; ^\ell\!\7h$ for a Cartan
subalgebra as in \Ruf{ZS} then each $\E_\jmath$ is spanned by some of
the vectors in the basis $(^{{\ell}}f_{j})$. For each fixed $\7g\cong
\su(p,q)$ there are only finitely many $G^0$-conjugacy classes of
centralizers $C_{\7g}(\7v\red)$ as $\7v\subset \7g$ varies through all
qualifying MASAs in $\7g$.\Formend

For the sake of clarity let us mention that in general there are
infinitely many conjugacy classes of qualifying MASAs $\7v$ while the
above Lemma asserts that there are only finitely many conjugacy
classes of the corresponding toral parts $\7v\red$. The point here is
that for a fixed $\7v\red$ there may be infinitely many non-conjugate
qualifying MANSAs in $C_{\7g}(\7v\red)\ss\,$.

\medskip The above lemma describes the structure of $\7v\red$ and its
centralizer $C_{\7g}(\7v\red)$ in an elemen\-tary-geometric way and
shows that both determine each other uniquely. For the description of
$\7v=\7v\red\oplus\7v\nil$ it is therefore enough to determine all
possible $\7v\nil$. These split into a direct sum 
$$\eqalign{\7v\nil&=\bigoplus_{\5K\cup\5L}\7n_{\jmath}\Steil{with}\cr
\7n_{\jmath}:&=\7v\nil \cap\7s_{j}\steil{for}\7s_{\jmath}:=\cases{\;\;
\su(\2E_{\jmath},\6h_{\jmath})& $\jmath\in\5K$\cr \noalign{\vskip4pt}
\big(\7{sl}(\2E_\jmath)\oplus \7{sl}(\2E_{\bar\jmath})
\big)^\sigma&$\jmath\in\5L$ .\cr }\cr}\Leqno{VN}$$ Every
$\7n_{\jmath}$ is an abelian $\ad$-nilpotent subalgebra of
$\7s_{\jmath}$. In case $\jmath\in\5K$ the algebra $\7n_{\jmath}$ has
dimension $p_{\jmath}+q_{\jmath}-1$.  As a consequence, $p_{\jmath}=0$
is possible only if $q_{\jmath}=1$ (since in case $q_{\jmath}>1$ the
form $\6h_{\jmath}$ is definite and every $\ad$-nilpotent element in
$\7s_{\jmath}$ is zero). In the same way $q_{\jmath}=0$ implies
$p_{\jmath}=1$. In case $\jmath\in\5L$ the algebra $\7n_{\jmath}$ has
dimension $2m_{\jmath}-2$.

\medskip
Before we turn to the corresponding classification result we need to
extract some invariants from the equations in \ruf{EI}. For every set
$A$ denote by $\5F(A)$ the {\sl free commutative monoid} over $A$. We
write the elements of $\5F(A)$ in the form $\sum_{\alpha\in
A}n_{\alpha}\!\cd\alpha$ with $n_{\alpha}\in\NN$ and
$\sum_{A}n_{\alpha}<\infty$. Here we use the free monoids over the
following sets, where $\NN=\{0,1,2,.\dots\}$:
$$\eqalign{\3K:&=\{(s,t)\in\NN^{2}:(st=0)\Rightarrow(s+t=1)\}\,,\cr
\3L:&=\NN\backslash\{0\}\Steil{and}\3J:=\3K\cup
\3L\,.\cr}\Leqno{MO}$$ Then\hfill
$\5D:=\5F(\3J)=\5F(\3K)+\5F(\3L)\,$,\hfill\phantom{Then}\hfill

\noindent and the permutation of $\3J$ defined by $(s,t)\mapsto(t,s)$
on $\3K$ and the identity on $\3L$ induces an involution
$\3D\mapsto\3D\opp$ of $\5D$. As an example, the {\sl opposite} of
$\3D=4\cd(3,5)+2\cd\97$ is $\3D\opp=4\cd(5,3)+2\cd\97$. Notice that
$2\cd\97$ and $7\cd\92$ are different elements in
$\5F(\3L)\subset\5D$. The set $\3J$ can be considered in a canonical
way as subset of $\5F(\3J)$ by identifying $\3j\in \3J$ with $1\cd
\3j\in\5D$, but for better distinction we write $\3j$ instead of $1\cd
\3j$ only if no confusion is likely. Also, for better distinction we
write the natural numbers in $\NN\backslash\{0\}$ in boldface if we
consider them as element of $\3L$.

\Definition{}
For every $p,q\in\NN$ we denote by $\5D_{p,q}\subset\5D$ the subset of
all
$$\3D=\sum_{\3j\in\3J}n_{\3j}\cd\3j\;\in\;\5F(\3J)\qquad\hbox{satisfying}$$ 
$$p=\sum_{\3j=(s,t)\in \3K}n_{\3j}s+\sum_{\3j\in
\3L}n_{\3j}j\Steil{and}q=\sum_{\3j=(s,t)\in \3K}n_{\3j}t+\sum_{\3j\in
\3L}n_{\3j}j\,,$$ where $j$ is the natural number underlying
$\3j$. Then $\5D_{p,q}\opp=\5D_{q,p}$ and
$\5D_{p,q}+\5D_{p',q'}\subset\5D_{p+p',q+q'}$ are obvious.

\medskip To every qualifying MASA $\7v\subset\su(p,q)$ we associate an
element $\3D(\7v)$ of $\5D$ that only depends on the
$\SU(p,q)$-conjugation class of $\7v$ and is called the {\sl
$D$-invariant of $\7v$}: Suppose that $\7v$ (after a suitable
conjugation) gives rise to the equations \Ruf{ZR} and \Ruf{EI}. Then
just put
$$\3D(\7v):=\sum_{\jmath\in\5K}1\cd(p_{\jmath},q_{\jmath})+
\sum_{\jmath\in\5L}1\cd \3m_{\jmath}\,.\Leqno{DF}$$ For instance, the
CSAs in $\7g=\su(p,q)$ are precisely the qualifying MASAs
$\7v\subset\7g$ with $\3D(\7v)\in\5F(A)$ with
$A:=\{(1,0),(0,1),\91\}$. Indeed, in the notation of \Ruf{ZS} we have
$\3D(\,{}^{\ell}\!\7h)=(p-\ell)\cd(1,0)+(q-\ell)\cd(0,1)+\ell\cd\91$.

The relevance of the $D$-invariants
for our classification problem is demonstrated by the following two
results.  Recall that $G=N_{\SL(p{+}q,\CC)}(\su(\E,\6h_{p,q}))$ and $
G^0=\SU(\2E,\6h_{p,q})$ with $G\ne G^0$ only if $p=q$.

\Proposition{RP} $\3D\in\,\5D$ is the $D$-invariant of a qualifying
MASA $\7v$ in $\su(p,q)$ if and only if $\3D\in\5D_{p,q}$.\Formend

\Proposition{NC} Let $\7v_{1},\7v_{2}$ be two qualifying MASAs in
$\7g=\su(p,q)$. Then the toral parts $\7v\red_{1},\7v\red_{2}$ (and
hence also the corresponding centralizers) are $G^{0}$-conjugate
in $\7g$ if and only if $D(\7v_{1})=D(\7v_{2})$. In case $p=q$ the
toral parts $\7v\red_{1},\7v\red_{2}$ are $G$-conjugate if and only if
$D(\7v_{1})=D(\7v_{2})$ or $D(\7v_{1})=D(\7v_{2})\opp$ (in this case
$\SU(\2E,\6h_{p,q})$ has index two in $G$).\Formend

Our classification problem now reduces to the following task:
For every $p,q\ge1$ with $p+q\ge3$ and every $\3D$ in the finite set
$\5D_{p,q}$ determine all $G$-conjugacy classes of qualifying MASAs
$\7v\subset\su(p,q)$ with $\3D(\7v)=\3D$.

\Joker{}{Explicit classification for small values of $q$:}
For $\5D_{p,q}$ in the cases $p\ge q=1,2$ we have the
following explicit lists (without repetitions).

\noindent $\5D_{p,1}$ consists of all invariants\0
$(p-s)\cd(1,0)+(s,1)$ for $s=1,2,\dots,p$. \1
$(p-1)\cd(1,0)+\91\Steil{and}p\cd(1,0)+(0,1)$.

\smallskip\noindent $\5D_{p,2}$ consists of all invariants \1
$\91+\5D_{p-1,1}$, \1 $(p-s-t)\cd(1,0)+(s,1)+(t,1)$ for all $1\le s\le
t$ with $s+t\le p\,$,\1 $(p-s)\cd(1,0)+(s,2)\,$,
$\;(p-s)\cd(1,0)+(0,1)+(s,1)$ for $1\le s\le p\,$,\1
$(p-2)\cd(1,0)+\92\,$, $\;p\cd(1,0)+2\cd(0,1)\,$.

\noindent Notice that in $\5D_{2,2}$ there exists an invariant that is not
self-opposite, e.g. $(1,0)+(1,2)$.  As a consequence, in case
$\7g\cong\su(2,2)$ there are eleven $\SU(2,2)$-conjugation classes of
centralizers in contrast to the only ten $G$-conju\-gation classes
in this case.

\Joker{}{Nilpotent parts of qualifying MASAs.}  So far we have given a
description of conjugacy classes of the toral parts $\7v\red$ of
qualifying MASAs in $\7g$. For the description of $\7v=\7v\red\oplus
\7v\nil$ it is therefore sufficient to determine all possible
nilpotent parts $\7v\nil$. According to \ruf{OL}, each such $\7v\nil$
is a maximal abelian nilpotent subalgebra (MANSA) of
$C_{\7g}\ss(\7v\red)$, compare \Ruf{MD}. On the other hand, given any
qualifying MASA $\7v\subset \7g$, each MANSA $\7n$ of $
C_{{\7g}}\ss(\7v\red)$ in the $-1$-eigenspace of $\tau$ gives rise to
a qualifying MASA $\7v_{\7n}:=\7v\red\oplus \7n$ in $\7g$. Given one
of the finitely may conjugacy classes of centralizers $C_{\3D}$ with
$\3D\in \5D_{p,q}$, our task is therefore to classify the MANSAs in
the semisimple part $C_\3D\ss$ of $C_\3D$. As already explained,
$C_\3D\ss$ decomposes uniquely into simple ideals $\7g_\jmath$, each
of them being isomorphic either to $\su(p_\jmath,q_\jmath)$ or to
${\7{sl}}_{m_\jmath}(\CC)$. Consequently each qualifying MANSA
$\7n\subset C_\3D\ss$ has the unique decomposition
$\7n=\bigoplus_{\5J_{1}}\!\7n_{\jmath}$ with
$\5J_{1}:=\{\jmath\in\5K\cup\5L:\dim\E_{\jmath}>1\}\subset\5J$,
compare \ruf{EI}.  Each of the factors $\7n_{\jmath}$ is a qualifying
MANSA in $\7g_{\jmath}$, more precisely:
$$\diagram{ \llap{$C_{\7g}(\7v\red) = $ }C_{\3D}&=& \7v\red &\oplus &
\bigoplus_{\5J_{1}}& \7g_{\jmath}\cr \cup & & \Vert\;\; & & & \cup\cr
\7v_{\7n}&=& \7v_{\7n}\red &\oplus & \bigoplus_{\5J_{1}}& \7n_{\jmath}
} \qquad\qquad
\diagram{\7g_\jmath&\;\cong\su(p_\jmath,q_\jmath)\hfill\cr
\noalign{\vskip4pt} &\hbox{or}\hfill\cr \noalign{\vskip4pt}
\7g_\jmath&\;\cong \5R^\CC_\RR({\7{sl}}(m_{\jmath},\CC))\;.}\kern-3em
\Leqno{ZF}$$

\bigskip Summarizing the results of the present subsection, our next
task is to determine maximal abelian $\ad$-nilpotent subalgebras
$\7n_{\jmath}\subset\7g_{\jmath}$ with either
$\7g_{\jmath}=\su(\E_\jmath,\6h_\jmath)\cong
\su(p_{\jmath},q_{\jmath})$ or $\7g_\jmath=({\7{sl}}(\E_\jmath)\oplus
{\7{sl}}(\E_{\overline\jmath}))^\sigma\cong
\5R^\CC_\RR({\7{sl}}_{m_\jmath}(\CC))$. Here we can restrict to
subalgebras that are contained in the $(-1)$-eigenspace of an
involution $\tau$ coming from a conjugation on the vector spaces
$\E_\jmath$ and $\E_\jmath\oplus \E_{\overline\jmath}$ respectively.
In the following we discuss the cases
${\7{sl}}(\E)^\sigma=\su(\2E,\6h)\cong\su(p,q)$ and
$\big(\7{sl}(\2E)\oplus
\7{sl}(\overline\2E)\big)^\sigma\cong{\7{sl}}(m,\CC)$ separately.

\KAP{MANS}{MANSAs in $\su(p,q)^{-\tau}$}

For notational simplicity let us drop the subscript $\jmath\,$ for the
rest of this section and write $\E=\E_j$ as well as $(p,q)=
(p_\jmath,q_\jmath)$ for the type of the restriction of $\6h$ to
$\E_\jmath$ and also $\7g=\su(\E,h)$. As before we denote by $\sigma$
the conjugation on the complexification $\7g^{\CC}=\7{sl}(\E)$ with
$\Fix(\sigma)=\7g$.  For every $z\in\End(\E)$ we denote by
$z^{\star}\in\End(\E)$ the adjoint with respect to $\6h$.

Without loss of generality we assume $pq\ne 0$, since otherwise
$p+q=1$ and thus $\su(E,\6h)=0$, see \Ruf{MO}.  As before, $m=p+q$ and
$r=m-1$.

\medskip
In order to classify the maximal nilpotent Lie subalgebras 
$\7n=\7v\nil\subset \su(E,\6h)$ we relate them to 
nilpotent commutative and {\sl associative} $\RR$-algebras, 
see the Appendix for the terminology. 

\smallskip\Proposition{QA} Let $\7n\subset \7g=\su(\2E,\6h)$ be a Lie
subalgebra. Then the following conditions are equivalent: \0 $\7n$ is
maximal among all abelian Lie subalgebras $\7a$ of $\7g$ such that
every element of $\7a$ is a nilpotent endomorphism of $\2E$.  \1 $\7n$
is maximal among all abelian subalgebras of $\7g$ that are
ad-nilpotent in $\7g$.  \1 The complexification $\7n^{\!\CC}$ is
$\sigma$-stable and is maximal among all abelian and nilpotent
associative subalgebras of $\End(\E)$.

\noindent If these conditions are satisfied then
$\7n^{\!\CC}\!\oplus\; \CC\id_{\2E}$ is a maximal abelian subalgebra
of the complex {\rm associative} algebra $\End(\2E)$. Further, in the
notation of \ruf{QD}$\,$ff. for every $\7n\subset \7g$ satisfying one
(and hence all) of the conditions {\rm(i) - (iii)} the following
holds: \1 The subspaces $\2K_{\7n^\CC}$ and $\2B_{\7n^\CC}$ are
$\6h$-orthogonal in $\E$, implying $d_1=d_3$.  \1 There exists an
$\7n^{\!\CC}\!$-adapted decomposition $\E_1\oplus \E_2\oplus \E_3$ of
$\E$ such that $\E_1$ and $\E_3=\2K$ are $\6h$-isotropic and
$\E_2\perp_h(\E_1\oplus \E_3)$. If, in addition, $\7n$ is
$\tau$-stable then the adapted decomposition can be chosen to be
$\tau$-stable, too.  The restriction of $\6h$ to ${\E_2}$ is of type
$(p-d_1, q-d_1)$.  \Formend

\medskip\noindent Now the conjugation $\tau$ of $\E$ comes into
play. For $\2V:=\E^{\tau}\cong\RR^{m}$ we have $\E=\2V\oplus i\2V$,
and we identify $\End(\2V)$ in the obvious way with the real
subalgebra $\End(\2E)^{\tau }$ of $\End(\2E)$. A crucial observation
is the following refinement of Proposition \ruf{QF}.

\Proposition{FQ} Assume that $\7n\subset \su(\E,\6h)^{-\tau}$ is
maximal among abelian and ad-nilpotent subalgebras of $\su(E,\6h)$. As
before, let $\5N:=i\7n\subset\End(\E)$. Then, with the notation of
Proposition \ruf{QF}, the following holds:\0 $\dim \Ann(\5N)=1$. In
particular, $\dim \2V_1=\dim \2V_3=1$ for any $\5N\!$-adapted
decomposition of $\,\2V$.  \1 Fix generators $v_1\in \2V_1$ and
$v_3\in \2V_3$ with $\6h(v_1,v_3)=1$. This yields canonical
identifications $\2V_1=\RR=\2V_3$, $\,\Ann(\5N)=\RR$,
$\,\5N_{21}=\Hom(\RR,\2V_2)=\2V_2$ and
$\5N_{32}=\Hom(\2V_2,\RR)=\2V_2^{*}$ (the dual of $\2V_{2}$). The map
$\6J:\2V_{2}\to\2V_{2}^{*}$ is given by $\6J(y)(x)=\6h(x,y)$ for all
$x,y\in\2V_{2}$. With all these identifications the matrix
presentation in Proposition \ruf{QF} reads
$$\5N=\left\{\pmatrix{0&0&0\cr y&\6N(y)&0\cr t&\6J(y)&0\cr
}:y\in\2V_{2},\,t\in\RR\right\}\;\subset
\S(\2V,\6h)\subset\End(\2E)\;,$$ where $ \S(\2V,\6h) \subset\End(\2V)$
is the linear subspace of all $\6h$-selfadjoint operators on $\2V$.
\1 The restriction of $\6h$ to $\2V_2$ has type $(p-1, q-1)$ and
$N(y)\in \S(\2V_{2},\6h)$ for every $y\in\2V_{2}$.

\Proof (ii): From \ruf{QA} and \ruf{QF}.(ii) follows that for a
maximal abelian and ad-nilpotent subalgebra $\7n\subset \su(\E,\6h)$
we have $\Ann(i\7n)=\{x\in \Ann(\7n^{\!\CC})=\Hom(\E_1,\E_3):
x=x^{\star}\}$. Since at the same time $\7n$ is contained in the
$(-1)$-eigenspace of $\tau$ the first part of the lemma together with
\ruf{QF}.(ii) imply
$\Ann(i\7n)=\Hom(\2V_1,\2V_3)=\Hom(\E_1,\E_3)^\tau$. This is only
possible if $\dim \E_1=\dim \E_3=1$. \qed

\medskip The next proposition shows that the classification of maximal
nilpotent subalgebras $\7n\subset \su(\E,\6h)$, contained also in
${\7{sl}}(\E)^{-\tau}$, reduces to the classification of abstract
associative nilpotent subalgebras (of $\End(\2E)$) with 1-dimensional
annihilator. Crucial for the following theorem is the construction of
a non-degenerate 2-form $\6b=\6b_{\pi}$ depending on a suitable
projection $\pi$, see \ruf{BN}. Keeping also in mind Proposition
\ruf{QE} and Lemma \ruf{UG} we have:

\Theorem{ZG} Let $\5N$ be an arbitrary commutative associative and
nilpotent $\RR$-algebra with $\dim \Ann(\5N)=1$ and let $\2V:=\5N^{0}$
be its unital extension.  Fix an identification $\Ann(\5N)=\RR$ and a
projection $\pi$ on $\2V$ with range $\Ann(\5N)=\RR$ satisfying
$\pi(\One)=0$. Then for the left regular representation $L$ of
$\2V=\5N^{0}$ and the symmetric real 2-form $\6b:\2V\times\2V\to\RR$ we
have:\0 $L(\5N)$ is a maximal nilpotent and abelian subalgebra of
$\End(\2V)$ contained in $\S(\2V,\6b)$. \1 Let $\E:=\2V\oplus i\2V$ be
the complexification of $\,\2V$ and $\tau$ the conjugation on $\2E$
with $\2E^{\tau}=\2V$. Furthermore, denote the unique hermitian
extension of $\6b$ to $\2E\times\2E$ by the same symbol $\6b$.
Then
$\7n:=iL(\5N)$ is a subset of $ \su(\E,\6b)^{-\tau}$ and is
a maximal abelian and ad-nilpotent Lie subalgebra of
$\su(\E,\6b)$. Finally, $\exp\7n^\CC\subset\SL(\E)$ has an
open orbit in $\PP(\E)$. \1 Every maximal abelian and ad-nilpotent
subalgebra of $\su(\E,\6h)$ which is also contained in
${\7{sl}}(\E)^{-\tau}$ for $\tau$ as in \ruf{SE} is equivalent to some
$iL(\5N)$ as above.  \Formend

\bigskip\medskip\noindent Suppose that for $k=1,2$ there are given two
abelian nilpotent associative $\RR$-algebras $\5N_{k}$ each with
1-dimensional annihilator $\5A_{k}$ and assume that the corresponding
2-forms $\6b_{k}$ on the corresponding unital extensions have types
$(p_{k},q_{k})$ with respect to the linear isomorphisms
$\lambda_{k}:\5A_{k}\cong\RR$. Then for
$\5I:=\{(x,y)\in\5A_{1}\oplus\5A_{2}:
\lambda_{1}(x)=\lambda_{2}(y)\}$ the quotient algebra
$\5N:=\qu{(\5N_{1}\oplus\5N_{2})}{\5I}$ is an abelian nilpotent
associative algebra with 1-dimensional annihilator
$\5A:=\qu{(\5A_{1}\oplus\5A_{2})}{\5I}$, and the 2-form $\6b$ on $\5N$
induced by $\6b_{1}\times\6b_{2}$ on $\5N_{1}\times\5N_{2}$ has type
$(p_{1}+p_{2}-1,q_{1}+q_{2}-1)$.

\bigskip\Joker{}{MANSAs in $\su(p,q)$ for low values of $q$.} By our
above considerations, to determine all MANSAs $\7n\subset\su(p,q)$ up
to $\SU(p,q)$-conjugacy it is equivalent to determine up to
isomorphism all real nilpotent abelian associative algebras $\5N$ with
annihilator $\5A$ such that for some linear isomorphism $\5A\cong\RR$
the form $\6b$ has type $(p,q)$ on $\5N^{0}$.  For low values of $q$
this can be done:

\kap{$q=1$:} There is precisely 1 equivalence class of MANSAs
in $\su(p,1)^{-\tau}$ for every $p\ge1$. Indeed, for every $\5N$ with
annihilator $\5A\cong\RR$ the factor algebra $\qu{\5N}{\5A}$ must be
a zero product algebra.

\kap{$q=2$:} There are precisely $\min(p,3)$ equivalence classes of
MANSAs in $\su(p,2)^{-\tau}$ for every $p\ge1$. Representing algebras
$\5N$ are obtained as follows: For every $n\ge1$ with $n\le\min(p,3)$
let $\5N_{1}$ be the cyclic abelian algebra of dimension $n$, compare
Example \ruf{EX}, and identify $t\in\RR$ with $t\xi^{n}$ in the
annihilator $\5A_{1}$ of $\5N_{1}$. Then the corresponding form
$\6b_{1}$ has type $(1,1)$, $(2,1)$, $(2,2)$ for $n=1,2,3$
respectively. Next choose an abelian nilpotent algebra $\5N_{2}$ with
1-dimensional annihilator $\5A_{2}$ such that the construction
$\5N:=\5N_{1}\oplus\5N_{2}/\5I$ as above leads to a nilpotent algebra
with 1-dimensional annihilator $\5A$ such that the corresponding
2-form $\6b$ on $\5N^{0}$ has type $(p,2)$. This is always possible
since $\qu{\5N_{2}}{\5A_{2}}$ must be a zero product algebra.

\KAP{MANSA}{MANSAs in $\7{sl}(m,\CC)$}

In this section we deal with the simple factors
$\7g_{\!\jmath}\cong{\7{sl}}(m_\jmath,\CC)$ in \Ruf{EI}. We retain our
convention from the last section and drop the index $\jmath$ from our
notation, that is, we consider abelian nilpotent subalgebras of
$\End(\2E)$ that are contained in $\7g=\7{sl}(m,\CC)$, where the
latter space is considered as a real Lie algebra.  Recall that in this
case the restriction of the involution $\tau$ to $\7g$ is given by the
map $x\mapsto -x'$, where $x'$ is the adjoint with respect to the
complex bilinear non-degenerate symmetric 2-form $\beta$ given by
$\beta(v,w):=\6h(v,\tau w)$ for all $v,w\in\2E$, compare \Ruf{NT}.
Note that the complexification $\7g^{\!\CC}$ is isomorphic to the
product ${\7{sl}}(m,\CC)\times {\7{sl}}(m,\CC)$.
 
\smallskip\Proposition{QB} Let $\7g=\5R^\CC_\RR({\7{sl}}(\E))$ and let
$\7n\subset\7g$ be a Lie subalgebra.  Then the following conditions
{\rm (i) -- (iii)} are equivalent: \0 $\7n$ is maximal among
ad-nilpotent and abelian $\RR$-subalgebras of ${\7{sl}}(\E)$.  \1
$\7n$ is maximal among all abelian $\CC$-subalgebras of ${\7{sl}}(\E)$
that are ad-nilpotent in ${\7{sl}}(\E)$.  \1 $\7n$ is maximal among
abelian and nilpotent subalgebras of the associative complex algebra
$\End(\E)$.

\noindent If these conditions are satisfied, then $\7n^{\!\CC}\oplus\;
\CC\cd\id_{\2E}$ is a maximal abelian subalgebra of the complex
associative algebra $\End(\2E)$. Furthermore we have: \1 Let $\7n$ be
maximal among ad-nilpotent and abelian $\CC$-subalgebras of
${\7{sl}}(\E)$ contained in ${\7{sl}}(\E)^{-\tau}$.  Then
$\dim\Ann(\7n)=1$, i.e., $1=\dim \E_1=\dim \E_3$ for any $\7n$-adapted
decomposition of $\2E$ {\rm(dimensions over $\CC$)}.  \1
Identifying $\CC=\E_1=\E_3=\Ann(\7n)$, the matrix presentation in
Proposition \ruf{QF} reads
$$\5N=\left\{\pmatrix{0&0&0\cr y&\6N(y)&0\cr t&\6J(y)&0\cr
}:y\in\2E_{2},\,t\in\CC\right\}\;\subset
\S(\2E,\beta)\subset\End(\2E)\;,$$ where $ \S(\2E,\beta)
\subset\End(\2E)$ is the linear subspace of all $\beta$-selfadjoint
operators on $\2E$.

\Proof (iv): Proposition \ruf{QF} implies that $\Ann
(\7n)=\Hom(\E_1,\E_3)$ for a given $\7n$-adapted decomposition of
$\E$. On the other hand, $\7n$ is also contained in
${\7{sl}}(\E)^{-\tau}$, i.e., $x=x'$ for all $x\in {\7{sl}}(\E)$.
Hence, we can choose an $\7n$-adapted decomposition $\E=\E_1\oplus
\E_2\oplus \E_3$ such that $\E_1$ and $\E_3$ are
$\beta$-isotropic and $(\E_1\oplus \E_3)$, $\E_2$ are
$\beta$-orthogonal.  In the matrix presentation of $\7n$ as in
Proposition \ruf{QF} the two above conditions imply
$\Ann(\7n)=\Hom(\E_1,\E_3)=\{x\in \Hom(\E_1,\E_3): x=x'\}$. This is
only possible if $\dim E_1=\dim \E_3=1=\dim \Ann(\7n)$.\qed

\medskip Similar to Theorem \ruf{ZG} we have

\Theorem{TS} Let $\5N$ be an arbitrary associative $\CC$-algebra which
is commutative nilpotent and has annihilator $\Ann(\5N)$ of dimension
1. Let furthermore $\pi$ be an arbitrary (complex linear) projection
on its unital extension $\2E:=\5N^{0}$ with range $\Ann(\5N)$ and
$\pi(\One)=0$, and fix an identification $\Ann(\5N)=\CC$.  Then
$L(\5N)$ is maximal in the class of all nilpotent and abelian
subalgebras $A\subset\End(\E)$ which are contained in
$\S(\E,\6b)_{\tr=0}={\7{sl}}(\E)^{-\tau}$. Here $\6b=\6b_{\pi}$ is as
in \ruf{BN}, $x'$ is the adjoint with respect to the complex bilinear
2-form $\6b$ and $\tau:\End (\E)\to \End(\E)$ is given by $x\mapsto
-x'$.  On the other hand, every maximal abelian and ad-nilpotent
subalgebra of ${\7{sl}}(\E)$ which is also contained in
${\7{sl}}(\E)^{-\tau}$ for $\tau$ as in \ruf{NT} is equivalent to some
$L(\5N)$ as above.  \Formend

\medskip\noindent Note that the complex nilpotent group $N^\CC$,
corresponding to $\7n^\CC\cong \7n\times \7n\subset {\7{sl}}(\E)\times
{\7{sl}}(\E)\cong ({\7{sl}}(\E))^\CC$ does not have an open orbit in
$\PP(\E\oplus \E)$, but the subgroup corresponding to the following
subalgebra does
$$(\7n\times \7n) \oplus \CC(\id,-\id)=(\7v\nil)^\CC\oplus
(\7v\red)^\CC\subset C_{{\7{sl}}(\E\oplus
\E)}(\7v\red)\subset{\7{sl}}(\E\oplus \E)\,.$$
\kap{MANSAs in $\7{sl(m,\CC)}$ for low values of $m$.} There
exist exactly 1,1,1,2,3 equivalence classes of qualifying MANSAs in
$\7{sl(m,\CC)}$ for $m=1,2, 3,4,5$ (see also the more detailed
description at the end of the following Section \ruf{Normal}).  With
our construction of nilpotent algebras out of cubic forms $\6c$ in
case $\FF=\CC$ (compare Proposition \ruf{HU}) it follows that there are
infinitely many equivalence classes of qualifying MANSAs in
$\7{sl(m,\CC)}$ for every $m\ge8$.

\KAP{Normal}{Normal forms for equations}

Every local tube realization $T_{F}=V+iF\subset E:=V\oplus iV$ of
$S_{p,q}$ is characterized by a qualifying MASA $\7v\subset\su(p,q)$.
In addition, the base $F$ of the tube can always be chosen to be a
closed (real-analytic) hypersurface in the real vector space $V$, see
\Lit{FKAP}. In the following we want to find canonical real-analytic
real valued functions $\psi$ on $V$ with $d\psi\ne0$ everywhere such
that $F=\{x\in V:\psi(x)=\psi(0)\}^{0}$, where the upper index
${}^{0}$ means to take the connected component containing the origin.
For this we consider the $D$-invariant $\3D(\7v)$ of
$\7v=\7v\red\oplus \7v\nil$, see \Ruf{DF}, and start with the special
case that $\3D(\7v)\in\3J=\3K\cup\3L$. The general case with
$\3D(\7v)\in\5F(\3J)$ arbitrary then is obtained by putting these
special equations together.

\kap{1. Case $\3D(\7v)\in\3K$:} Let $\3j:=(p,q)$ and
$n:=p+q-1$. Furthermore let $\hat V_{\3j}:=\RR^{p+q}$ with coordinates
$(x_{0},x_{1},\dots,x_{n})$ and define the linear form $\lambda_{\3j}$
on $\hat V_{\3j}$ by $\lambda_{\3j}(x)=(p+q)x_{0}$. Also, identify
$V_{\3j}:=\RR^{n}$ with coordinates $(x_{1},\dots,x_{n})$ in the
obvious way with the hyperplane $\{x\in\hat
V_{\3j}:\lambda_{\3j}(x)=0\}$. We define $\Psi_{\3j}$ as the set of
all real-analytic functions $\psi(x)=e^{x_{0}}f(x_{1},\dots,x_{n})$ on
$\hat V_{\3j}$ where $f$ is an extended real nil-polynomial on
$V_{\3j}$ and the second derivative of $\psi$ at the origin of $\hat
V_{\3j}$ has type $(p,q)$, compare the Appendix for the notion of a
nil-polynomial.

It is clear that for every real $t>0$ and every $g\in\GL(V_{\3j})$
with $e^{x_{0}}f(x)$ also the function $t\,e^{x_{0}}f(g(x))$ is
contained in $\Psi_{\3j}$. Furthermore, $\Psi_{\3j\opp}=-\Psi_{\3j}$
is evident for the opposite $\3j\opp=(q,p)$.  In particular,
$\Psi_{(1,0)}=\{te^{x_{0}}:t>0\}$,
$\Psi_{(1,1)}=\{te^{x_{0}}x_{1}:t\ne0\}$ and $\Psi_{(2,1)}$ is the
orbit of $\psi=e^{x_{0}}(x_{2}+x_{1}^{2})$ under the group
$\RR^{*}\times\GL(2,\RR)$.

\medskip \kap{2. Case $\3D(\7v)\in\3L$:} Let $\3j:=\3m$ for some
integer $m\ge1$ and put $n:=m-1$. Consider $\hat V_{\3m}:=\CC^{m}$
with complex coordinates $(z_{0},z_{1},\dots,z_{n})$ as real vector
space and define the linear form $\lambda_{\3m}\,$ on $\hat V_{\3m}\,$
by $\lambda_{\3m}(z):=m(z_{0}+\overline z_{0})$. Furthermore, consider
$W_{\3m}:=\CC^{n}$ with coordinates $(z_{1},\dots,z_{n})$ in the
obvious way as linear subspace of $\hat V_{\3m}$. We define
$\Psi_{\3m}$ as the set of all real-analytic functions
$\psi(z)=\Re\big(e^{z_{0}}f(z_{1},\dots,z_{n})\big)$ on $\hat V_{\3j}$
where $f$ is an extended complex nil-polynomial on $W_{\3m}$.  Then
the second derivative of every $\psi\in\Psi_{\3m}$ at the origin has
type $(m,m)$.

The group $\CC^{*}\!\times\GL(W_{\3m})$ acts in a canonical way on the
extended complex nil-polynomials on $W_{\3m}$ and thus also on
$\Psi_{\3m}$. In particular, $\Psi_{\3m}$ is the orbit of the
functions $\Re(e^{z_{0}})$, $\Re(e^{z_{0}}z_{1})$ and
$\Re(e^{z_{0}}(z_{2}+z_{1}^{2}))$ for $m=1,2,3$ respectively.

\medskip \kap{3. Case $\3D(\7v)$ arbitrary:} Then
$\3D:=\3D(\7v)\in\5D_{p,q}$ for integers $p,q\ge1$ with $p+q\ge3$, and
there exists a unique sum representation $\3D=\sum_{\alpha\in
A}\3j_{\alpha}$ with $(\3j_{\alpha})_{\alpha\in A}$ a finite family in
$\3J=\3K\cup\3L$.  Put $\hat V_{\D}:=\bigoplus_{\alpha\in A} \hat
V_{\3j_{\alpha}}$ and define the linear form $\lambda_{\D}$ on $\hat
V_{\D}$ by
$(x_{\3j_{\alpha}})\mapsto\sum_{\alpha}\lambda_{\3j_{\alpha}}(x_{\3j_{\alpha}})$.
Furthermore, let $\Psi_{\D}$ be the space of all functions
$$\psi:\;(x_{\3j_{\alpha}})_{\alpha}\;\longmapsto\;\sum_{\alpha}
\psi_{\alpha}(x_{\3j_{\alpha}})$$ on $\hat V_{\D}\,,$ where
$(\psi_{\alpha})_{\alpha\in A}$ is an arbitrary family of functions
$\psi_{\alpha}\in\Psi_{\3j_{\alpha}}$. The second derivative of every
$\psi\in\Psi_{\D}$ at the origin then has type $(p,q)$.

\bigskip The relevance of the vector spaces $\hat V_{\D}$ with linear
form $\lambda_{\D}$ and function space $\Psi_{\D}$ is the following:
Consider in $\hat V_{\D}$ the hyperplane $V:=\{x\in \hat
V_{\D}:\lambda_{\D}(x)=0\}$. Then for every $\psi\in\Psi_{\D}$ the
analytic hypersurface
$$F:=\big\{x\in V:\psi(x)=\psi(0)\big\}^{0}\Leqno{UE}$$ is the base of
a local tube realization of $S_{p,q}$ with $D$-invariant $\3D$, and
every local tube realization of $S_{p,q}$ with $D$-invariant $\3D$
occurs this way up to affine equivalence. Indeed, every local tube
realization of $S_{p,q}$ is associated with a qualifying MASA
$\7v\subset\su(p,q)$. In particular, for $\2E:=\CC^{p+q}$ the
complexification $\7v^{\CC}:=\7v\oplus\,i\!\7v\subset\7{sl}(\2E)$ has
an open orbit $\O$ in the projective space $\PP(\2E)$ that is the
image of the locally biholomorphic map $\phi:\7v^{\CC}\to\PP(\2E)$,
$\;\xi\mapsto\exp(\xi)a$, where $a$ is a suitable point in $\PP(\2E)$.
Then every connected component $M$ of $\phi^{-1}(S_{p,q})$ is a closed
tube submanifold of $\7v^{\CC}$, lets take the one that contains the
origin. Then the base $F:=M\cap i\7v$ of the tube manifold $M$ has
\Ruf{UE} as defining equation if we put $V:=i\7v\subset\7{sl}(\2E)$
and $\psi:=\phi_{|V}$. Now consider the extended space $\hat
V:=\RR\cd\id\,\oplus\;i\!\7v\subset\gl(\2E)$ and let $\tr$ be the
trace functional on $\hat V$. Also extend $\psi$ to $\hat V$ by
$t\cd\id\oplus\,x\mapsto e^{t}\psi(x)$ and denote the extension by the
same symbol $\psi$. 
\nline Now the $D$-invariant $\3D(\7v)$ and the
corresponding decomposition of $\7v\red$ in Lemma \ruf{EI} comes into
play.  We can identify $\hat V$ with $\hat V_{\D}$ and $\tr$ with
$\lambda_{\D}$ in a canonical way. Since $\psi$ is defined in terms of
$\exp$ it is compatible with the decomposition in Lemma \ruf{EI} and we
only have to discuss defining equations for the case that there is
only one summand in \Ruf{ZR}, that is, that $\3D(\7v)\in\3J$:

\noindent{\sl 1. Case $\3D(\7v)\in\3K$:} Let $\3D(\7v):=(p,q)$.  In
case $pq=0$ we have $V=0$ and we have up to a positive factor
$\psi(x)=(p-q)e^{x}$ on $\hat V=\RR$.  We therefore assume $p,q\ge1$
in the following. But then $\5N:=V=i\7v$ is an associative commutative
nilpotent real subalgebra of $\End(\2E)$ with 1-dimensional
annihilator $\5A$, compare the Appendix. Choose a pointing $\omega$ on
$\5N$ such that the associated symmetric 2-form $\6h(x,y)=\omega(xy)$
has type $(p,q)$ on $\5N^{0}$. Then as base for a tube realization
associated with $\7v$ we can take
$$F=\{x\in\5N:f(x)=0\}\Steil{with}f(x):=\6h(\exp x/2,\exp x/2)=
\omega(\exp x)\,.$$ But $f$ is an extended real nil-polynomial on
$\5N$ and $F$ is a smooth algebraic hypersurface.

\noindent{\sl 2. Case $\3D(\7v)\in\3L$:} Let $\3D(\7v):=\3m$ and put
$n:=m-1$.  Then
$\7v=\RR\cd\id+\,\5N\;\subset\;\gl(m,\CC)\;\subset\;\7u(m,m)$, where
$\5N\subset\7{sl}(m,\CC)$ is a complex MANSA and at the same time a
commutative associative nilpotent complex subalgebra of
$\End(\CC^{m})$. In case $m=1$ we have $\5N=0$ and
$\psi(z)=\Re(e^{z})$ on $\hat V=\CC$. Let us therefore assume $m>1$ in
the following. Then $\5N$ has annihilator $\5A$ of complex dimension
1. Let $\omega$ be a pointing on $\5N$. Then the real symmetric 2-form
$\6h(z,w)=\Re\,\omega(zw)$ has type $(m,m)$ on the complex unital
extension $\5N^{0}=\CC\cd\id+\,\5N\subset\End(\CC^{m})$.  We have to
consider the complexification $\7v^{\CC}=\7v\oplus\,i\!\7v$ and to
restrict the exponential mapping to $i\!\7v$. For this, we may
identify the $\RR$-linear space $V:=i\!\7v$ with
$i\RR\cd\id+\,\5N\;\subset\;\5N^{0}$ and get as base for a tube
realization associated to $\7v$ the hypersurface $$F=\{w\in
i\RR\cd\id+\,\5N:\6h(\exp w/2,\exp w/2)=0\}\,.$$ Writing
$w=is\cd\id+z$ with $s\in\RR$, $z\in\5N$ we have $$\6h(\exp w/2,\exp
w/2)=\Re(e^{is}f(z))\Steil{for}f(z):=\omega(\exp z)\,,$$ that is, $f$
is an extended complex nil-polynomial on $\5N\cong\CC^{n}$, and $F$ is
affinely equivalent to the non-algebraic hypersurface
$$\{(s,z)\in\RR\oplus\CC^{n}:\Re\big(e^{is}
f(z_{1},\dots,z_{n})\big)=0\}^{0}\,.$$
\vskip-28pt\qed

\bigskip \Joker{}{Local tube realizations corresponding to Cartan
subalgebras of $\su(p,q)$:} By the above we know $\Psi_{\3j}$ for all
$\3j\in \{(1,0),(0,1),\91\}$ and thus we can explicitly write down the
normal form equations of every CSA in $\su(p,q)$: For fixed $p,q\ge1$
and every $\ell\ge0$ with $\ell\le\min(p,q)$ consider the CSA
${}^{\ell}\!\7h$ as defined in Section \ruf{centralizers}. Then we
have

\smallskip\centerline{$\3D:=
\3D({}^{\ell}\!\7h)=(p-\ell)\cd(1,0)+(q-\ell)\cd(0,1)+\ell\cd\91$.}
\noindent With $d:=p+q-2\ell$ then
$$V_{\D}=\big\{(z,t)\in\CC^{\ell}\oplus\RR^{d}:
\sum_{k=1}^{\ell}(z_{k}+\overline z_{k})+\sum_{k=1}^{d}t_{k}=0\big\}$$
and as tube base we can take a connected component of the set of all
$(z,t)\in V_{\D}$ satisfying
$$\sum_{k=1}^{\ell}\Re(e^{z_{k}})\;+\;\sum_{k=1}^{p-\ell}
e^{t_{k}}\,=\;\sum_{k=1}^{q-\ell}e^{t_{p-\ell+k}}\,.\Leqno{OQ}$$

\medskip \Joker{HF}{Comparing with the equations of Isaev-Mishchenko:}

It is easy to write down explicitly $\Psi_{\3j}$ for all $\3j=(p,q)$
with $q\le2$, getting back the classifications in \Lit{DAYA} and
\Lit{ISMI}. In the following we compare the equations obtained in
\Lit{ISMI} with ours, where $\3D$ is the corresponding $D$-invariant
and $n=p\,$:

\noindent{\bf types 1\Kl, 4\Kl, 5\Kl:} $\,\3D=s\cd(1,0)+(n-s,2)$. The
corresponding MANSA in $\su(n-s,2)$ has nil-index 2,3,4 respectively.

\noindent{\bf type 2\Kl:} $\,\3D=s\cd(1,0)+(0,1)+(n-s,1)$.

\noindent{\bf type 3\Kl:} $\,\3D=s\cd(1,0)+(n-s-1,1)+\91$,

\noindent{\bf type 6\Kl:} $\,\3D=(n-2)\cd(1,0)+\92$. 

\noindent{\bf type 7\Kl:} $\,\3D=s\cd(1,0)+(t,1)+(n-s-t,1)$.

\noindent{\bf types 8\Kl, 9\Kl, 10\Kl:} These types correspond to the
three Cartan subalgebras of $\su(n,2)$ and are affinely equivalent to
the equations \Ruf{OQ} for $\ell=0,2,1$ respectively.

\KAP{Some}{Some Examples}
We will give applications of Proposition \ruf{HU} in the real as well
as in the complex case. We start with the real version.

\Joker{real}{Examples obtained from real cubic forms} Let $W$ be a
real vector space of dimension $2n$ and $\6q$ a 
quadratic form of type $(n,n)$ on $W$. Then there exists a
decomposition $W=W'\oplus W''$ into totally isotropic linear
subspaces. Let furthermore $\6c$ be a cubic form on $W'$ and 
define the function $f$ on $V:=W\oplus\RR$ by
$$f(x,y,t):=t+\6q(x+y)+\6c(x)\Steil{for all}t\in\RR,x\in W',y\in
W''\,.\Leqno{KH}$$ Then $f$ is an extended nil-polynomial on $V$ and
the hypersurface $F:=\{v\in V:f(v)=0\}$ in $V$ is the base of a local
tube realization for $S_{p,p}$ with $p=n+1$. On the other hand, $F$ is
affinely homogeneous and also the complement $V\backslash F$ is
affinely homogeneous, compare the end of the Appendix. The complement
$V\backslash F$ decomposes into two affinely homogeneous domains
$D^{\pm}$, the tube domains over these domains are complex affinely
homogeneous domains in $V^{\CC}=V\oplus iV$.  With \Ruf{HX} we see
that there exists a real $n\choose3$-parameter family of cubic forms
on $W'$ leading to pairwise affinely inequivalent examples (notice
that property $(*)$ of Proposition \ruf{HW} is satisfied for all
$\alpha(t_{j})$ in \Ruf{HX} near $\6c_{0}$). In particular this shows
that there are infinitely many affinely non-equivalent local tube
realizations for $S_{4,4}$.

\Joker{complex}{Examples obtained from complex cubic forms} We start
with a more general situation: Suppose that $V$ is a complex vector
space and $f:V\to\CC$ is a holomorphic submersion with $f(0)=0$. Then
$$F:=\{z\in V:f(z)=0\}^{0}$$ is a complex hypersurface in $V$, the
complement $D:=V\backslash F$ is a domain in $V$ and
$$F:=\big\{(t,z)\in\RR\oplus V:\Re\big(e^{it}f(z)\big)=
0\big\}^{0}\Leqno{KI}$$ is a real hypersurface in $\RR\oplus V$.  With
respect to the canonical projection $\pr:F\to V$, $(t,z)\mapsto z$,
the surface $F$ is a covering over the domain $D$ and a trivial real
line bundle over $H$.  In fact, $\tilde F:=\pr^{-1}(F)$ is a connected
real hypersurface in $F$, while the open subset $\tilde
D:=\pr^{-1}(D)$ in $F$ in general is not connected.  \nline Now assume
that $V=W\oplus\CC$ and that $f$ is an extended complex nil-polynomial
of degree $\le3$ on $V$ as considered in Proposition \ruf{HU}. As a
consequence of Proposition \ruf{KA} the group
$\6A:=\{g\in\Aff(V):f\circ g=f\}$ acts transitively on every level set
$f^{-1}(c)$ in $V$. Also, for every $s\in\CC$ the linear
transformation
$\theta_{s}:=e^{s}\id_{W'}\oplus\,e^{2s}\id_{W''}\oplus\,e^{3s}\id_{\CC}$
satisfies $f\circ\theta_{s}=e^{3s}f$, see also \Ruf{KB}.  The group
$\CC\times\6A$ acts by the affine transformations
$$(t,z)\mapsto\big(t-3\,\Im(s)\,,\,\theta_{s}\,g(z)\big),\qquad
s\in\CC,g\in\6A,$$ on $F$ and has precisely three orbits there -- the
closed orbit $\tilde F$ and the two connected components $\tilde
D^{\pm}$ of the domain $\tilde D$. Also, the translation
$(t,z)\mapsto(t+\pi,z)$ interchanges these two domains $\tilde D^{+}$,
$\tilde D^{-}$ in $F$.  \nline Putting things together we got the
following: Let $W$ be a complex vector space of dimension $2n$ and
$\6q$ a non-degenerate quadratic form on $W$. Then there exists a
decomposition $W=W'\oplus W''$ into totally isotropic linear
subspaces. Let furthermore $\6c$ be a cubic form on $W'$ and define
the function $f$ on $V:=W\oplus\CC$ by \Ruf{KH} with $\RR$ replaced by
$\CC$. Then $f$ is an extended complex nil-polynomial and $F$ as defined in
\Ruf{KI} is the base of a local tube realization $M\subset
U^{\CC}=U\oplus iU$ for $S_{p,p}$, where $U:=\RR\oplus V$ and
$p=2(n+1)$. Furthermore, the real hypersurface $M$ in $U^{\CC}$ contains
an affinely homogeneous domain. For every $n\ge3$ we get a complex
$n\choose3$-parameter family of pairwise affinely inequivalent
examples.

\KAP{Appendix}{Appendix -- Nilpotent commutative algebras}

\smallskip \noindent{\sl In this Appendix, all occurring algebras are
either associative or Lie. Throughout, $\FF$ is an arbitrary
base field of characteristic zero. For every associative
algebra $A$, every $x\in A$ and every integer $k\ge1$ we put
$$x^{(k)}:={1\over
k!}x^{k}\Steil{and}x^{(0)}:=\One\steil{if}A\steil{has a
unit}\One\,.\Leqno{KP}$$ }

\medskip\noindent We collect several purely algebraic statements that
are used in the paper and might be of independent interest. Some of
them are probably known to the experts. Since we could not find a
reference in the literature we state it here. Recall e.g. our convention
that $\End(\E)$ is the {\sl associative} endomorphism algebra while
$\gl(\E)$ is the same space, but endowed with the corresponding Lie
product. We start with a standard definition.

\Definition{GX} Let $\5N$ be an commutative associative algebra over
$\FF$ and define the ideals $\5N^{k}\subset\5N$ inductively by
$\5N^{1}=\5N$ and $\5N^{k+1}=\langle\5N\5N^{k}\rangle$. Then $\5N$ is
called {\sl nilpotent} if $\5N^{k+1}=0$ for some $k\ge0$, and the
minimal $k$ with this property is called the {\sl nil-index} of
$\5N$. Furthermore, $\5A:=\Ann(\5N):=\{x\in \5N: x\5N=0\}$ is called
the {\sl annihilator} of $\5N$.

\medskip\centerline{\bf The general embedded case} In the following,
let $\E$ be a vector space of finite dimension $m\ge2$ over $\FF$.  For
every subalgebra $\5N\subset\End(\E)$ define the following
characteristic subspaces of $\E$:
$$\2B:=\2B_{\5N}:=\langle \5N(v):v\in \E\rangle\Steil{and}
\2K:=\2K_{\5N}:=\{v\in \E: \5N(v)=0\}\,.\Leqno{QD}$$

\Proposition{RF} Suppose that $\5N$ is maximal among all commutative
and nilpotent subalgebras of $\;\End(\E)$. Then \0 $0\ne
\2K_{\5N}\subset \2B_{\5N}\ne \E$. Also, $\2K_{\5N}=\2B_{\5N}$ holds
if and only if $\5A=\Ann(\5N)$ has nil-index 1.  \1 $\5N\oplus \FF\cd
\id$ is maximal among all commutative subalgebras of $\,\End(\E)$.  \1
$\5N$ is irreducible on $\E$, i.e., for every $\5N$-invariant
decomposition $\2E=\2E'\oplus \2E''$ either $\2E'=0\,$ or
$\,\2E''=0$.\Formend

\Joker{}{$\5N$-adapted decompositions and matrix presentations.}  For
every $\5N$ satisfying the assumptions in Proposition \ruf{RF}
we select subspaces $\2E_1, \2E_2$ of $\E$ such that $\E_{1}\oplus
\2B_{\5N}=\E$ and $\2E_{2}\oplus \2K_{\5N}=\2B_{\5N}$.  Then, for
$\E_{3}:=\2K_{\5N}$, we have the decomposition
$$\E=\E_{1}\oplus \E_{2}\oplus \E_{3}\Steil{with}d_{j}:=\dim\2E_{j}\steil{for}j=1,2,3\,,\Leqno{ZC}$$ that we also call an
{\sl $\5N$-adapted decomposition}. Every $x\in\End(\E)$ can
be written as $3\times3$-matrix $(x_{jk})$ with
$x_{jk}\in\Hom(\E_{k},\E_{j})$.  With
$\pi_{jk}:\End(\E)\to\Hom(\E_{k},\E_{j})$ we denote the projection
$x\mapsto x_{jk}$.

We call two subalgebras $\5N\subset \End(\E)$ and $\5N'\subset \End
(\E')$ {\sl conjugate} if there exists an invertible $\Psi\in
\Hom(\E,\E')$ such that $\5N'=\Psi\circ \5N\circ\Psi^{-1}$.  One of
our goals is to decide under which conditions two isomorphic
subalgebras (isomorphic as abstract $\FF$-algebras) are already
conjugate. In general, there exist isomorphic subalgebras which are
not conjugate.

It is obvious that a nilpotent subalgebra $\5N\subset\End(\E)$
contains only nilpotent endomorphisms.  By a theorem of Engel the
converse is also true: A subalgebra $\5N\subset \End(\E)$ consisting of
nilpotent endomorphisms only is nilpotent and there is a full flag
$$F_{1}\subset F_2\subset\cdots \subset
F_m=\E,\qquad\dim F_{k}=k,$$ which is stable under $\5N$, i.e.,
with respect to a suitable basis of $\E$ the algebra $\5N$ consists of
strictly lower-triangular matrices in $\FF^{m\times m}$.  With this
notation we can state

\Proposition{QF} Suppose that $\5N\subset\End(\2E)$ satisfies the
assumptions of \ruf{RF} and has nil-index $\nu$. For a fixed
$\5N$-adapted decomposition $\E=\E_{1}\oplus \E_{2}\oplus \E_{3}$ and
all $1\le j,k\le3$ put $\5N_{jk}:=\pi_{jk}(\5N)$. Then:\smallskip \0
There exists a linear bijection $\6J:\5N_{21}\to\5N_{32}$ and a linear
map $\6N:\5N_{21}\to\5N_{22}$ such that
$$\5N=\left\{\pmatrix{0&0&0\cr y&\6N(y)&0\cr t&\6J(y)&0\cr
}:y\in\5N_{21},\,t\in\Hom(\2E_{1},\2E_{3})\right\}\;.$$\1
$\5A:=\{x\in\5N:x_{21}=0\}\cong\Hom(\2E_{1},\2E_{3})$ is the
annihilator of $\5N$.\1 $\5N_{21}\times\5N_{21}\to\5N_{31}$,
$\,(x,y)\mapsto \6J(x)\circ y$, is a non-degenerate symmetric 2-form
(in fact, is equivalent to the restriction of the form $\6b$ defined
in \Ruf{BN} after the obvious identifications).\1 $\5N_{22}$ is a
nilpotent commutative subalgebra of $\End(\2E_{2})$\nline with nil-index
$\le\max(\nu-2,1)$.  \1
$\2E_{2}=\langle\5N_{23}(\2E_{1})\rangle$ and
$\{z\in\2E_{2}:y(z)=0\steil{for all}y\in\5N_{32}\}=0$. \1
$d_{1}d_{3}+\lceil d_{2}/\mu\rceil\;\le\;\dim\5N\, \le [m^{2}/4]$ for
$m:=d_{1}+d_{2}+d_{3}=\dim\E$ and $\mu:=\min(d_{1},d_{3})$.  In
particular, if $d_{1}=1$ then $\dim \5N=d_{2}+d_{3}=m-1$. \Formend

\Remark{BD} The linear maps $\6N$ and $J$ in (i) above satisfy for all
$x,y\in\2E$ the relations: 
\item{(a)} $\6N(x)^{k}=0$ for some integer $k$, 
\item{(b)}$\6N(x)y=\6N(y)x$, \quad $\6J(x)y=\6J(y)x$ \ and \ $\6J(x)\6N(y)=\6J(y)\6N(x)$,
\item{(c)}$\6N(\6N(x)y)=\6N(x)\6N(y)$. 
\par\noindent
On the other hand, let a vector space $\2W$ over $\FF$ be given. Every
pair $\6N:\2W\to \End(\2W)$,\break
$\6J:\2W\buildrel\approx\over\to\2W^*$ of linear maps satisfying (a) -
(c) gives rise by (i) above to a commutative maximal nilpotent
subalgebra $\5N\subset\End(\2E)$ with $\2E=\FF\oplus \2W\oplus
\FF$, i.e., $\2E_1=\FF=\2E_3$, $\2E_2=\2W$ and
$\5N$ has 1-dimensional annihilator. In particular, $\6N\equiv0$ and
$\6J$ given by any symmetric and non-degenerate scalar product on $\2W$
trivially satisfy (a) - (c) and define a maximal nilpotent subalgebra $\5N$
with $\dim \5N=\dim \E-1$ and nil-index $2$. In this case $\5N/\5A$ is the zero product
algebra.

\Remark{} The upper bound in inequality (vi) is sharp as the nilpotent
subalgebra
$$\{x\in\End(\E):x_{jk}=0\steil{if}(j,k)\ne(3,1)\}$$ for
$d_{2}=d_{3}=\lceil m/2\rceil$ shows.  It is much harder to find
better lower bounds for $\dim \5N$, not to speak of sharp ones. For
infinitely many values of $m$ (starting with $m=14$)
there exist maximal commutative and
nilpotent subalgebras of $\End(\FF^{m})$ with $\dim \5N<m-1$, see \Lit{LAFF}.
\Formend

\bigskip \centerline{\bf Abstract commutative nilpotent algebras}
\medskip The left-regular representation of a nilpotent (associative)
algebra is not faithful, contrary to the case of any unital algebra.
For every nilpotent algebra $\5M$ denote by
$\5M^{0}:=\FF\cd\One\oplus\5M$ its unital extension. Such extensions
of nilpotent algebras are precisely those $\FF$-algebras, which
contain a maximal ideal of codimension 1 consisting of nilpotent
elements only. Denote by $L:\5M^{0}\to \End(\5M^{0})$ the
corresponding left regular representation. It is obvious that the
algebras $\5M$ and $L(\5M)\subset\End(\5M^{0})$ are isomorphic via
$L$.

\Proposition{QE} Let $\5M$ be a commutative nilpotent $\FF$-algebra of
 finite dimension. Then \0 $L(\5M)$ is maximal among all commutative
 nilpotent subalgebras of $\End(\5M^{0})$ and consists entirely of
 nilpotent endomorphisms. Furthermore, $\2K_{L(\5M)}=\Ann(\5M)$ and
 $\2B_{L(\5M)}=\5M$. The image of the unital extension, $L(\5M^{0})$ is
 maximal among all commutative subalgebras of $\End (\5M^{0})$.  \1 A
 maximal commutative nilpotent subalgebra $\5N\subset \End(\E)$ is
 conjugate to the image $L(\5M)$ of some commutative nilpotent algebra
 $\5M$ as above if and only if $\codim_{\E}\2B_{\5N}=1$, see \Ruf{QD}
 for the notation.  \1 Let $\5N\subset\End(\E)$ and $\5M\subset
 \End(\2F)$ be two maximal commutative nilpotent subalgebras. If
 $\codim_\E\2B_{\5N}=\codim_{\2F}\2B_{\5M}=1$ then $\5M$ and $\5N$ are
 conjugate (in the above defined sense) if and only if $\5M$ and $\5N$
 are isomorphic as abstract algebras.  There exist non-conjugate
 subalgebras $\5N,\5M\subset \End(\E)$ with $\codim_{\E}\2B _{\5N}>1$,
 which are isomorphic as abstract algebras.  \Formend

\Joker{}{Associated 2-forms.}  Let $\5N$ be an commutative and
nilpotent $\FF$-algebra and $\5A:=\Ann(\5N)$ its annihilator. On the
unital extension $\o{\5N}=\FF\cd\One\oplus\5N$ fix a projection
$\pi=\pi^{2}\in\End(\o{\5N})$ with range $\pi(\5N)=\5A$ and
$\pi(\One)=0$. Then $$ \6b_{\pi}:\5N^{0}\times\5N^{0}\to\5A,\qquad
(x,y)\longmapsto\pi(xy)\Leqno{BN}$$ defines an $\5A$-valued symmetric
2-form.  Clearly, the restriction of $\6b_{\pi}$ to $\5N$ factorizes
through $\5N/\5A\times\5N/\5A$ and we write also $\6b_{\pi}$ for the
corresponding $\5A$-valued 2-form on $\5N/\5A$.  Furthermore, $\pi$
determines the decomposition
$$\eqalign{ \5N^{0}&=\5N_1\oplus\5N_2\oplus\5N_3\Steil{with}\cr
\5N_{1}&=\FF\cd\One,\qquad \5N_{2}=\5N\cap \ker
\pi\cong\5N/\5A\Steil{and} \5N_{3}=\5A=\Ann(\5N)\;.\cr }\Leqno{ZD}$$
The canonical isomorphism $\cong$ in \Ruf{ZD} makes $\5N_{2}$ to an
algebra that we denote by $\5N_{2}^{\pi}$. In terms of $\pi$ and the
algebra structure on $\5N$ the product on $\5N_{2}^{\pi}$ is given by
$(x,y)\mapsto (\id-\pi)(xy)$ for $x,y\in\5N_{2}^{\pi}$.

\medskip

\Lemma{UG} Let $\5N\ne 0$ be a commutative nilpotent associative
$\FF$-algebra with annihilator $\5A:=\Ann(\5N)$ and let $L:\5N^{0}\into
\End(\5N^{0})$ be the left regular representation of the unital
extension $\5N^{0}$.  Fix a projection $\pi$ on $\5N^{0}$ with range $\5A$
as above and let $\6b_{\pi}$ be the associated 2-form \Ruf{BN}. Then \0 The
decomposition \Ruf{ZD} is an $L(\5N)$-adapted decomposition of
$\5N^{0}$.  \1 $\6b_{\pi}$ is nondegenerate.  \1 $\6b_{\pi}$ is associative, or
equivalently, every $L(y)\in\End(\5N^{0})$ is $\6b_{\pi}$-selfadjoint.  \1 The
subspaces $\5N_{1}$ and $\5N_{3}$ are $\6b_{\pi}$-isotropic while the
subspaces $(\5N_1\oplus \5N_{3})$ and $\5N_2$ are $\6b_{\pi}$-orthogonal to
each other. Consequently, the restriction of $\6b_{\pi}$ to the algebra
$\5N_{2}^{\pi}\cong\qu{\5N}{\5A}$ is a non-degenerate associative
$\5A$-valued 2-form. \1 In the particular case $\FF=\RR$ and
$\dim\5A=1$ the following holds after fixing a linear isomorphism
$\psi:\5A\to\RR$: The type $(p,q)$ of the real 2-form $\psi\circ\6b_{\pi}$
on $\5N^{0}$ does not depend on the choice of the projection
$\pi$. \Formend

\bigskip \medskip \centerline{\bf Commutative nilpotent algebras
consisting of selfadjoint endomorphisms}

\medskip\noindent As indicated in the main part of this paper the
classification of the various local tube realizations of $S_{p,q}$ is
equivalent to the classification of maximal abelian subalgebras $\7v$
in $\su(p,q)=\su(\E,\6h)$, contained in the $(-1)$-eigenspace of the
involution $\tau$ of the Lie algebra $\su(p,q)$, up to conjugation by
elements from $G:=N_{\SL(\E)}(\su(p,q))$, compare \ruf{TA}.

Every such abelian subalgebra $\7v$ admits a unique decomposition
into its ad-reduc\-tive and its ad-nilpotent part,
$\7v=\7v\!\red\oplus \7v\!\nil$.  Also, we have the further
decomposition $\7v\!\nil=\bigoplus_\jmath \7n_\jmath$, where the
building blocks $\7n_j$ are various ad-{\sl nilpotent} abelian
subalgebras $\7n_{\!\jmath}$ maximal in $\su(p_\jmath,q_\jmath)$ or
$\7{sl}(m_{\jmath},\CC)$. In turn, such Lie algebras (also contained
in $\Fix(-\tau)\,$) are in a 1-1-correspondence to associative
commutative and nilpotent subalgebras of $\End(\2V_\jmath)$, that
consist of selfadjoint endomorphisms with respect to a symmetric
non-degenerate 2-form $\6h$ on $\2V$ over $\FF=\RR$ or $\FF=\CC$.  This
is one motivation to investigate a symmetric version of maximal
commutative and nilpotent subalgebras in $\End(\E)$.

In this subsection let $\2E$ be an arbitrary vector space of finite
dimension over $\FF$ and $\6h:\2E\times \2E\to \FF$ a symmetric
non-degenerate 2-form. With $\S(\E,\6h)\subset \End(\E)$ we denote the
linear subspace of all operators $a$ that are selfadjoint with respect
to $\6h$ (that is, $\6h(ax,y)=\6h(x,ay)$ for all $x,y\in \2E$). Given
$\2V\subset \E$ we write $\2V^\perp$ for the orthogonal complement
with respect to $\6h$. Note that in general $\2V\cap \2V^\perp\ne
0$. For (maximal) nilpotent commutative subalgebras in $\End(\E)$,
contained in $\S(\E,\6h)$, there are $\5N$-adapted decompositions which
are in addition related to the 2-form $\6h$:

\Lemma{ND} Let $\6h:\E\times \E \to \FF$ be a non-degenerate symmetric
2-form and $\5N\subset \End(\E)$ a maximal nilpotent commutative
subalgebra, contained in $\S(\E,\6h)$. Let $\2K_\5N$, $\2B_\5N$ be the
characteristic subspaces as defined in \ruf{QD}.  Then: \0
$\2K^{}_{\5N}=\2B_{\5N}^{\perp}$.  \1 $\dim\Ann(\5N)=1$ 
\1 There exists an $\5N$-adapted
decomposition $\E=\E_1\oplus \E_2\oplus \E_3$ such that
\itemitem{\rm(a)} $\dim \E_1=\dim \E_3=1$. 
 \itemitem{\rm(b)} $\6h\big|_{\E_3}\!=0$,
$\6h\big|_{\E_1}\!=0$, and the pairing $\6h:\E_1\times \E_3\to \FF$ as
well as the restriction $\6h\big|_{\E_2}$ are non-degenerate.
\itemitem{\rm(c)} $\E_2=(\E_1\oplus \E_3)^{\perp}$.  \Formend

\Definition{} Let $\5N\subset \End(\E)$ be as in the previous Lemma.
We call every decomposition $\E=\E_1\oplus \E_2\oplus \E_3$ which
satisfies the condition (ii) in Lemma \ruf{ND} an $(\5N,\6h)$-{\it
adapted decomposition of} $\E$.

\medskip In the following let $\5N\subset \S(\E,\6h)$ be a maximal
nilpotent commutative subalgebra of $\End(\E)$ and $\5A:=\Ann(\5N)$ its
annihilator (which is 1-dimensional according to the above lemma).
Next, we relate the 2-form $b_\pi:\5N^{0} \times \5N^{0} \to
\5A:=\Ann(\5N)$ to the symmetric 2-form $\6h:\E\times \E\to \FF.$ Recall
that the choice of the projection $\pi:\5N\to \5A $ is equivalent to
the choice of a linear subspace $\5N_2\subset \5N$ with
$\5N=\5N_2\oplus \5A$. It is easy to see that every
$(\5N,\6h)$-adapted decomposition of $\E$ gives rise to the
complementary subspace $\5N_2:=\{n\in \5N: n(\E_1)\subset \E_2\}$,
i.e., $\5N_2\oplus \5A=\5N$. It turns out to be more subtle to prove
the opposite statement as it involves the solution of certain
quadratic equations in $\E$.

\Proposition{DX} Let $\5N\subset \S(\E,\6h)$ be a maximal nilpotent
commutative subalgebra. Then: \0 For every linear subspace
$\5N_2\subset \5N$ satisfying $\5N=\5N_2\oplus
\5A$, there exists an
$(\5N,\6h)$-adapted decomposition $\E=\E_1\oplus \E_2\oplus \E_3$
with
$$\5N_2=\{n\in \5N: n(\E_1)\subset \E_2\}\,.$$ \1 For $\E=\E_1\oplus
\E_2\oplus \E_3$ as in {\rm(i)} choose generators $e_1\in \E_1$,
$e_3\in \E_3$ with $\6h(e_1, e_3)=1$ and define $\kappa:\5A\to \FF$ by
$n(e_1)=\kappa(n)e_3$ for all $n\in \5A$.  Let $\pi:\5N^{0}\to \5A$
the projection corresponding with kernel $\5N_{1}\oplus\5N$ and
$\kappa\circ\6b_\pi:\5N^{0}\times \5N^{0}\to \FF$ the corresponding
nondegenerate symmetric 2-form.  Then the map
$$\5N^{0}\to \E,\quad m\mapsto m(e_1)$$ is an isometry between
$(\5N^{0}, \kappa\circ\6b_\pi)$ and $(\E,\6h)$ which respects the
decompositions $\5N^{0}=\FF\cd\One \oplus \5N_2\oplus\5A$ and
$\E=\E_1\oplus \E_2\oplus \E_3$\Formend

\bigskip\noindent Proposition \ruf{DX} is the main ingredient in the
proof of the following theorem, which can be considered as a symmetric
version of Proposition \ruf{QE}:

\Theorem{} Let $\6h$ and $\6h'$ be two non-degenerate symmetric forms on
$\E$ and $\E'$ respectively over $\FF$ , and let $\5N\subset
\End(\E)$, $\;\5N'\subset \End(\E')$ be two maximal nilpotent and
commutative subalgebras.  \0 Assume, that in addition
$\5N\subset\S(\E,\6h)$ and $\5N'\subset\S(\E',\6h')$ holds.  Then $\5N$
and $\5N'$ are isomorphic as $\FF$-algebras if and only if there
exists an isometry $\Psi:(\E,\6h)\to (\E',\6h')$ with $\5N'=\Psi\circ
\5N\circ\Psi^{-1}$.  \1 In particular, if $\,\E=\E'$ and $\6h=\6h'$, two
such subalgebras $\5N$, $\5N'$ contained in $\S(\E, \6h)$ are isomorphic
if and only if they are conjugate by an element in $\SO(\E,\6h)$.
\Formend

\medskip\noindent In case $\FF=\RR,\CC$ the above theorem has the
following application for the classification of maximal abelian
subalgebras of $\su(p,q)$ and $\7{sl}(m,\CC)$.  Let $\su(p,q)\cong
\su(\E,\6h)$ and $\tau:\su(p,q)\to \su(p,q)$ be as in \ruf{TA},
induced by a conjugation $\tau:\E\to \E$. Recall that we write
$\2V=\E^\tau$ for the real points with respect to $\tau$.  Note that
$\su(p,q)^\tau\cong \so(p,q)$ and $\SU(p,q)^\tau\cong \SO(p,q)\cong
\SO(\2V,\6h|_{\2V}).$ Further $\7{sl}(m,\CC)^\tau\cong \so(m,\CC)$,
i.e., $\tau: \7{sl}(m,\CC)\to \7{sl}(m,\CC)$ is induced by a symmetric
non-degenerate 2-form $\6h_\tau$ on $\E$.

\Corollary{} Two maximal abelian Lie subalgebras $\7v_1$,
$\7v_2$ in $\su(\E,\6h)$ respectively $\7{sl}(\E)$ consisting of
nilpotent elements and contained in $\su(\E,\6h)^{-\tau}\cong i\S(\2V,
\6h)_0$ respectively $\7{sl}(\E)^{-\tau}\cong \S(\E,\6h_\tau)_0$ are
conjugate under $\SO(\2V,\6h)$ respectively under
$\SO(\E,\6h_\tau)$ if and only if the corresponding associative
algebras $i\!\7v_1$, $i\!\7v_2$ in $\End(\2V)$ (respectively  $\7v_1$ and
$\7v_2$ in $\End(\E)$) are isomorphic as $\FF$-algebras.  \Formend

\bigskip\medskip\centerline{\Bf Some affinely homogeneous surfaces}

Let $\5N\ne0$ be a commutative associative nilpotent algebra over
$\FF$ of finite dimension with nil-index $\nu$. For every $k\ge1$
choose a linear subspace $V_{k}\subset\5N^{k}$ with
$\5N^{k}=V_{k}\oplus\5N^{k+1}$.  Then $\5N=\bigoplus_{k\ge1}V_{k}$
with $V_{\nu}=\5N^{\nu}$, and we write every $x\in\5N$ in the form
$x=\sum_{k=1}^{\nu}x_{k}$ with $x_{k}\in V_{k}$. With
$\pi:\5N\to\5N^{\nu}$ we denote the canonical projection $x\mapsto
x_{\nu}$. As before we denote by $\5N^{0}=\FF\cd\One\oplus\5N$ the
unital extension of $\5N$ and extend $\pi$ linearly to $\5N^{0}$ by
requiring $\pi(\One)=0$.

Denote by $\6P$ the space of all polynomial maps $\6p:\5N\to\5N$ of
the form $$x\;\longmapsto\;\sum \6p_{i_{1}i_{2}\dots
i_{r}}x_{i_{1}}x_{i_{2}}\cdots x_{i_{r}}\,,\Leqno{QI}$$ where the
integers $r\ge1$ and $\;1\le i_{1}\le\dots\le i_{r}$ satisfy
$\sum_{j=1}^{r}i_{j}\le\nu$ and the coefficients $\6p_{i_{1}i_{2}\dots
i_{r}}$ are from $\5N^{0}$. It is clear that with respect to composition
$\6P$ is a unital algebra over $\5N^{0}$. We are mainly interested in
polynomials $\6p\in\6P$ that are invertible in $\6P$, that is, where
for every $1\le j\le \nu$ the coefficient $\6p_{j}$ in front of the
linear monomial $x_{j}$ is invertible in $\5N^{0}$.

For every $\6p\in\6P$ the composition $\6f:=\pi\circ\6p$ is a
polynomial map $\5N\to\5N^{\nu}$ of degree $\le\nu$. Furthermore, in
case $\6p$ is invertible, the algebraic subvariety
$$F:=\{x\in\5N:\6f(x)=0\}\Leqno{RB}$$ is smooth, in fact, is the graph
of a polynomial map $\ker(\pi)\to\5N^{\nu}$. Denote by $\Aff(\5N)$ the
group of all affine automorphisms of $\5N$.

\Proposition{KA} Suppose that $\5N$ has nil-index $\nu\le4$ and that
$\6p\in\6P$ is invertible. In addition assume that \0 $\6p_{12}$ is
invertible in $\5N^{0}$ if $\nu=3$, \1 $V_{j}V_{k}\subset V_{j+k}$ for
all $j,k$ and that $\6p_{112}$, $\6p_{13}$ are invertible in $\5N^{0}$
if $\nu=4$.

\noindent Then for $\6f:=\pi\circ\6p$ the group
$\6A:=\{g\in\Aff(\5N):\6f\circ g=\6f\}$ acts transitively on every
translated subvariety $c+F=\6f^{-1}(c)$, $\;c\in\5N^{\nu}$.

\Proof We assume $\nu=4$, the cases $\nu<4$ are similar but easier.
For every $k=2,3,4$ denote by $\5L_{k}\subset\End(\5N)$ the subspace
of nilpotent transformations $x\mapsto \alpha_{1}x_{1}$, $\;x\mapsto
\alpha_{2}x_{1}+\alpha_{1}x_{2}$, $\;x\mapsto
\alpha_{3}x_{1}+\alpha_{2}x_{2}+\alpha_{1}x_{3}$ respectively with
arbitrary coefficients $\alpha_{j}\in V_{j}$.\nline Now fix an
arbitrary point $a\in F$ and denote by $\tau\in\Aff(\5N)$ the
translation $x\mapsto x+a$. A simple computation shows $$
\6f\circ\tau(x)=\6f(x)+x_{1}^{2}R_{2}(x)+ x_{1}R_{3}(x)+R_{4}(x)$$ for
suitable $R_{k}\in\5L_{k}$. Then
$\rho:=\id-\6p_{112}^{-1}R_{2}\in\GL(\5N)$ is unipotent and satisfies
$$ \6f\circ\tau\circ\rho\,(x)=\6f(x)+x_{1}S_{3}(x)+S_{4}(x)$$ for
suitable $S_{k}\in\5L_{k}$. Further
$\sigma:=\id-\6p_{13}^{-1}S_{3}\in\GL(\5N)$ satisfies
$$ \6f\circ\tau\circ\rho\circ\sigma(x)=\6f(x)+T_{4}(x)$$ for a
suitable $T_{4}\in\5L_{4}$. Finally,
$g(x)=\tau\circ\rho\circ\sigma(x-\6p_{4}^{-1}T_{4}(x))$ defines an
element $g\in\6A$ with $g(c)=c+a$ for all $c\in\5N^{\nu}$.\qed

\medskip \Remark{} The proof of Proposition \ruf{KA} also works for
fields $\FF$ of arbitrary characteristic. Special
polynomials $\6p\in\6P$ can be defined in the following way: Let
$$\Phi:=\sum_{k=1}^{\infty}c_{k}T^{k}\;\in\;\FF[[T]]\Leqno{OZ}$$ be an
arbitrary formal power series over $\FF$ with vanishing constant
term. Since $\5N$ is nilpotent $\6p(x):=\Phi(x)=\sum c_{k}x^{k}$ defines
a polynomial map $\6p\in\6P$. Clearly, $\6p$ is invertible if and only
if $c_{1}\ne0$. Furthermore, if we assume that $\FF$ has
characteristic $0$, then invertibility of $\6p_{12}$ is equivalent to
$c_{2}\ne0$ and invertibility of $\6p_{13}\6p_{112}$ is equivalent to
$c_{2}c_{3}\ne0$. For simplicity we may add a constant term $c_{0}$ to
the formal power series $\Phi$ (which will not count) if we at the same
time extend the projection $\pi$ from $\5N$ to its unital extension
$\5N^{0}=\FF\cd\One\oplus\5N$ by requiring $\pi(\One)=0$.  Later we are
mainly interested in the case where $\Phi=\exp$ is the usual exponential
series.

\medskip The conditions (i), (ii) in Proposition \ruf{KA} cannot be
omitted: As a simple example with $\nu=3$ consider the 3-dimensional
cyclic algebra $\5N$ with basis $e_{1}$, $e_{2}=e_{1}^{2}$,
$e_{3}=e_{1}^{3}$ satisfying $e_{1}^{4}=0$. Then, identifying $\5N$
with $\FF^{3}$ in the obvious way, we get for $\Phi=T+T^{3}$ in \Ruf{OZ}
that $\6f(x)=x_{3}+x_{1}^{3}$ on $\FF^{3}$. In this case, the group
$\6A$ does not act transitively on $F=\6f^{-1}(0)$ in general. In
fact, in case $\FF=\RR$ the affine group
$\Aff(F)=\{g\in\Aff(\5N):g(F)=F\}$ has precisely two orbits in $F$ -
the line $\RR e_{2}\subset F$ and its complement in $F$. Indeed, the
subgroup of all $(x_{1},x_{2},x_{3})\mapsto (tx_{1},x_{2}+s,t^{3}x_{3})$
with $s\in\RR$, $t\in\RR^{*}$ acts transitively on the complement.

\bigskip \medskip \centerline{\bf Nil-polynomials}
\medskip

In this subsection let $\5N$ be an arbitrary commutative associative
nilpotent algebra of finite dimension over $\FF$ with annihilator
$\5A$ of dimension $1$ and nil-index $\nu$. Let us call every linear
form $\omega$ on $\5N$ with $\omega(\5A)=\FF$ a {\sl pointing} of
$\5N$. Also, $\,\5N$ with a fixed pointing is called a {\sl pointed
algebra,} a PANA for short. Two PANAs $(\5N,\omega)$,
$(\tilde\5N,\tilde\omega)$ are called {\sl isomorphic} if there is an
algebra isomorphism $g:\tilde\5N\to\5N$ with $\tilde\omega=\omega\circ
g$. As before with projections we always consider every pointing
$\omega$ on $\5N$ linearly extended to $\5N^{0}$ by requiring
$\omega(\One)=0$.

For every vector space $V$ of finite dimension we denote by $\FF[V]$ the
algebra of all ($\FF$-valued) polynomials on $V$.

\Definition{RS} $f\in\FF[W]$ is called a {\sl nil-polynomial} on $W$
if there exists a PANA $(\5N,\omega)$ and a linear isomorphism
$\phi:W\to\ker(\omega)\subset\5N$ such that
$f=\omega\circ\exp\circ\,\phi$. Two nil-polynomials $f\in\FF[W]$,
$\tilde f\in\FF[\tilde W]$ are called {\sl equivalent} if there exists
$t\in\GL(\FF)\cong\FF^{*}$ and a linear isomorphism $g:W\to\tilde W$
with $\tilde f=t\circ f\circ g^{-1}$.

\Definition{HV} In case $V\ne0$ we call $f\in\FF[V]$ an {\sl extended
nil-polynomial} on $V$ if there exists a PANA $(\5N,\omega)$ and a
linear isomorphism $\phi:V\to\5N$ with
$f=\omega\circ\exp\circ\,\phi$. In case $V=0$ every constant in
$\FF^{*}$ is called an extended nil-polynomial on $V$.

\medskip
Nil-polynomials on vector spaces $W$ of dimension $n$ and extended
nil-polynomials on vector spaces $V$ of dimension $n+1$ correspond to
each other. Indeed, every extended nil-polynomial on $V$ is a sum
$f=\sum_{k>0}f_{[k]}$ of homogeneous parts $f_{[k]}$ of degree
$k$. Furthermore, $V=W\oplus A$ with $W=\ker f_{[1]}$ and $A=\{y\in
V:f_{[2]}(x+y)=0\;\forall x\in V\}\cong\FF$. Then the restriction of $f$
to $W$ is a nil-polynomial on $W$ and every nil-polynomial on $W$
occurs this way. For our applications in Section \ruf{Normal} we need
extended nil-polynomials. In the following we consider only
nil-polynomials for simplicity.

\medskip By definition, every equivalence class of nil-polynomials in
$\FF[W]$ is an orbit of the group $\GL(\FF)\times\GL(W)$ acting in the
obvious way on $\FF[W]$. For every pair of nil-polyno\-mials
$P\in\FF[W]$, $\tilde P\in\FF[\tilde W]$ we get a new nil-polynomial
$P\oplus\tilde P\in\FF[W\oplus \tilde W]$ by setting $(P\oplus\tilde
P)(x,\tilde x):=P(x)+\tilde P(\tilde x)$ for all $x\in W$ and
$x\tilde\in\tilde W$.

\medskip Fix a nil-polynomial $f\in\FF[W]$ in the following. Then we
have the expansion $f=\sum_{k\ge2}f_{[k]}$ into homogeneous
parts. Notice that $f_{[2]}$ is a non-degenerate quadratic form on
$W$. For every $k\ge2$ there is a unique symmetric $k$-form
$\omega_{k}$ on $W$ with
$$\omega_{k}(x,\dots,x)=k!\,f_{[k]}(x)\Leqno{HO}$$ for all $x\in
W$. Using $f_{[2]}$ and $f_{[3]}$ we define a commutative (not
necessarily associative) product $(x,y)\mapsto x\cd y$ on $W$ by
$$\omega_{2}(x\cd y,z)=\omega_{3}(x,y,z)\steil{for all}z\in
W\Leqno{HY}$$ and also a commutative product on $W\oplus\,\FF$ by
$$(x,s)(y,t):=(x\cd y,\omega_{2}(x,y))\,.\Leqno{HR}$$ Then, if
$f=\omega\circ\exp\circ\,\phi$ for a PANA $(\5N,\omega)$ with kernel
$\5K=\ker(\omega)$ and linear isomorphism $\phi:W\to\5K$ we have
$\omega_{k}(x_{1},\dots,x_{k})=\omega((\phi x_{1})(\phi
x_{2})\cdots(\phi x_{k}))$ for all $k\ge2$ and $x_{1},\dots,x_{k}\in
W$. For the annihilator $\5A$ of $\5N$ there is a unique linear
isomorphism $\psi:\FF\to\5A$ such that $\pi=\psi\circ\omega$ is the
canonical projection $\5K\oplus\5A\to\5A$. With these ingredients we
have

\Lemma{HZ} With respect to the product \Ruf{HR} the linear map
$$W\oplus\,\FF\;\to\;\5N\,, \quad(x,s)\mapsto\phi(x)+\psi(s)\,,$$ is
an isomorphism of algebras. In particular, $W$ with product $x\cd y$
is isomorphic to the nilpotent algebra $\5N/\5A$ and has nil-index
$\nu{-}1$. 

\Proof For all $x,y\in W$ we have

\smallskip\centerline{$(\phi(x)+\psi(s))(\phi(y)+\psi(t))=(N-A)+A$ with}

\smallskip\centerline{$N:=\phi(x)\phi(y)\in\5N\steil{and}A:=
\pi(\phi(x)\phi(y))=\psi(\omega_{2}(x,y))\in\5A$.}

\noindent It remains to show $N-A=\phi(x\cd y)$. But this follows from

\smallskip\centerline{$N-A\,\in\,\5K\Steil{and}\omega\big(\phi(x\cd
y)\phi(z)\big)=\omega\big(\phi(x)\phi(y)\phi(z)\big)=
\omega\big((N-A)\phi(z)\big)$} \noindent for all $z\in W$. \qed

\Corollary{HP} Every nil-polynomial $f$ on $W$ is uniquely determined
by its quadratic and cubic term, $f_{[2]}$ and $f_{[3]}$. In fact, the
other $f_{[k]}$ are recursively determined by
$$\omega_{k+1}(x_{0},x_{1},\dots,x_{k})=\omega_{k}(x_{0}\cd
x_{1},x_{2},\dots,x_{k})\Leqno{HQ}$$ for all $k\ge2$ and all
$x_{0},\dots, x_{k}\in W$, where the symmetric $\omega_{k}$ are
determined by \Ruf{HO}.\Formend 
Another application of \Ruf{HZ} is the following

\def\cdt{{\buildrel{\lower5pt\hbox{$\scriptscriptstyle\sim$}}\over\cd}}

\Proposition{} Let $(\5N,\omega)$, $(\tilde\5N,\tilde\omega)$ be
arbitrary PANAs and let $f\in\FF[W]$, $\tilde f\in\FF[\tilde W]$ be
associated nil-polynomials respectively. Then, if $f$ and $\tilde f$ are
equivalent as nil-polynomials, also $\5N$ and $\tilde\5N$ are isomorphic
as algebras.

\Proof Write $\tilde f=t\circ f\circ g^{-1}$ as in Definition \ruf{RS}
and define the products $\cd$ and $\cdt$ on $W$ and $\tilde W$ as in
\Ruf{HY}. With respect to these products $g:W\to\tilde W$ is an
algebra isomorphism. As in \Ruf{HR} the products $\cd$ and $\cdt$ extend to the algebras $W\oplus\FF\cong\5N$ and $\tilde W\oplus\FF\cong\tilde\5N$.
Finally $g\oplus\id$ gives an algebra isomorphism between them.\qed

\Lemma{} For every nil-polynomial $f$ on $W$ the cubic term $\6c:=f_{[3]}$
is trace-free with respect to the quadratic term $\6q:=f_{[2]}$, see 
{\rm \Lit{EAEZ} \p20} for this notion of trace,

\Proof Let $f\in\FF[W]$ be given by the PANA $(\5N,\omega)$ with
nil-index $\nu$. Without loss of generality we assume that
$W=\ker(\omega)$. Choose a basis $e_{1},\dots,e_{n}$ of $W$ and a
mapping $\alpha:\{1,\dots,n\}\to\{1,\dots,\nu{-}1\}$ such that
$\{e_{i}:\alpha(i)=\ell\}$ is a basis of $\5N^{\ell}/\5N^{\ell+1}$ for
$\ell=1,\dots,\nu{-}1$.  With respect to this basis the forms
$\6q$, $\6c$ are given by the tensors $g_{ij}=\omega_{2}(e_{i},e_{j})$
and $h_{ijk}=\omega_{3}(e_{i},e_{j},e_{k})$. It is clear that
$g_{ij}=0$ holds if $\alpha(i)+\alpha(j)\;>\;\nu$. Since $\alpha$
is surjective, this implies $g^{ij}=0$ if
$\alpha(i)+\alpha(j)\;<\;\nu$, where $(g^{ij})=(g_{ij})^{-1}$.  On
the other hand, $h_{ijk}=0$ if $\alpha(i)+\alpha(j)\;\ge\;\nu$,
proving the claim.\qed

 Corollary \ruf{HP} suggests the following question: Given a
non-degenerate quadratic form $\6q$ and a cubic form $\6c$ on $W$.
When does there exist a nil-polynomial $f\in\FF[W]$ with $f_{[2]}=\6q$
and $f_{[3]}=\6c\,$? Using $\6q$, $\6c$ we can define as above for
$k=2,3$ the symmetric $k$-linear form $\omega_{k}$ on $W$ and with it
the commutative product $x\cd y$ on $W$. A necessary and sufficient
condition for a positive answer is that $W$ with this product is a
nilpotent and associative algebra. As a consequence we get for every
fixed non-degenerate quadratic form $\6q$ on $W$ the following
structural information on the space of all nil-polynomials $f$ on $W$
with $f_{[2]}=\6q$: Denote by $\6C$ the set of all cubic forms on
$W$. Then $\6C$ is a linear space of dimension ${n+2\choose3}$,
$n=\dim W$, and
$$\6C_{\6q}:=\{\6c\in\6C:\exists\steil{nil-polynomial $f$ on $W$
with}f_{[2]}=\6q,f_{[3]}=\6c\}\Leqno{QH}$$ is an algebraic subset. The
orthogonal group $\O(\6q)=\{g\in\GL(W):\6q\circ g=\6q\}$ acts from the
right on $\6C_{\6q}$. The $\O(\6q)$-orbits in $\6C_{\6q}$ are in
1-1-correspondence to the equivalence classes of nil-polynomials $f$
on $W$ with $f_{[2]}=\6q$.

Examples of nil-polynomials of degree $3$ can be constructed in the
following way.

\Proposition{HU} Let $W$ be an $\FF$-vector space of finite dimension
and $\6q$ a non-degenerate quadratic form on $W$. Suppose furthermore
that $W=W_{1}\oplus W_{2}$ for totally isotropic linear subspaces
$W_{k}$ and that $\6c$ is a cubic form on $W_{1}$. Then, if we extend
$\6c$ to $W$ by $\6c(x+y)=\6c(x)$ for all $x\in W_{1}$, $y\in W_{2}$,
the sum $f:=\6q+\6c$ is a nil-polynomial on $W$. In particular,
$\6c\in\6C_{\6q}$.\nline Every $g\in\GL(W_{1})$ extends to an
$h\in\O(\6q)\subset\GL(W)$ in such a way that with
$\tilde\6c:=\6c\circ g$ also $\6q+\tilde\6c=f\circ h$ is a
nil-polynomial on $W$.

\Proof $\omega_{3}(x,y,t)=0$ for all $t\in W_{2}$ implies $W_{1}\!\cd
W_{1}\subset W_{2}$ and $W_{1}\!\cd W_{2}=0$, that is, $(x\cd y)\cd
z=0$ for all $x,y,z\in W$.\nline Now fix $g\in\GL(W_{1})$.  There
exists a unique $g^{\sharp}\in\GL(W_{2})$ with
$\omega_{2}(gx,y)=\omega_{2}(x,g^{\sharp}y)$ for all $x\in W_{1}$ and
$y\in W_{2}$. But then $h:=g\times (g^{\sharp})^{-1}\in\O(\6q)$ does
the job.\qed

\bigskip Let $\FF\subset\KK$ be a field extension and consider every
polynomial $f$ on $V$ in the canonical way as polynomial $\tilde f$ on
$V\otimes_{\FF}\KK$.  Then with $f$ also $\tilde f$ is a
nil-polynomial. In general, for non-equivalent nil-polynomials $f$,
$g$ on $V$ the nil-polynomials $\tilde f$, $\tilde g$ may be
equivalent. We use these extensions in case $\RR\subset\CC$.

\bigskip \medskip \centerline{\bf Graded PANAs}

Let $(\5N,\omega)$ with annihilator $\5A$ be a PANA in the
following. A grading then is a decomposition
$$\5N=\bigoplus_{k>0}\5N_{k}\,,\qquad \5N_{j}\5N_{k}\subset
\5N_{j+k}\;.\Leqno{HG}$$ Clearly $\5A=\5N_{d}$ for
$d:=\max\{k:\5N_{k}\ne0\}$ and $\nu=d/m$ is the nil-index of
$\5N$, where $m:=\min\{k:\5N_{k}\ne0\}$. Without loss of generality
we assume that $W:=\bigoplus_{k<d}\5N_{k}$ is the kernel of $\omega$.

For the corresponding nil-polynomial $f=\omega\circ\exp\;\in\;\FF[W]$
and $\ell:=d-1$ we then have
$$f=\sum_{k=2}^{\nu}f_{[k]}\Steil{with}f_{[k]}(x)={1\over
k!}\Big(\sum_{j_{1}+\dots+j_{k}=d}x_{j_{1}}x_{j_{2}}\cdots
x_{j_{k}}\Big)\Leqno{HI}$$ for all
$x=(x_{1},\dots,x_{\ell})\in\5N_{1}\oplus\cdots\oplus\5N_{\ell}=W\subset\5N$,
where every index $j_{\ell}$ in \Ruf{HI} is positive.  If we put,
using \Ruf{KP},
$$\Vert\mu\Vert
=\mu_{1}+2\mu_{2}+\dots+\ell\mu^{}_{\ell}\quad\Steil{and}\quad
x^{(\mu)} =x_{1}^{(\mu_{1})}x_{2}^{(\mu_{2})}\cdots
x_{\ell}^{(\mu_{\ell})}
$$ for every multi-index $\mu\in\NN^{\ell}$ and every
$x=(x_{1},\dots,x_{\ell})$, we can rewrite \Ruf{HI} as
$$f(x)=\sum_{\Vert\mu\Vert=d}x^{(\mu)}\,.\Leqno{QL}$$ 
We consider an example.

\Joker{EX}{Cyclic PANAs} For fixed integer $\nu\ge1$ let $\5N$ be the
{\sl cyclic} algebra of nil-index $\nu$ over $\FF$, that is, there is
an element $\xi\in\5N$ such that the powers $\xi^{k}$, $1\le k\le
\nu$, form a basis of $\5N$ and $\xi^{\nu+1}=0$. Then $\5N$ is a
graded algebra with respect to $\5N_{k}:=\FF\,\xi^{k}$ for all $k>0$
in \Ruf{HG} and becomes a PANA with respect to the pointing $\omega$
uniquely determined by $\omega(\xi^{k})=\delta_{k,\nu}$ for all
$k$. For $n:=\nu-1$ we identify $\FF^{n}$ and $\ker(\omega)$ via
$(x_{1},\dots,x_{n})=\sum x_{k}\xi^{k}$. The corresponding
nil-polynomial $f\in\FF[x_{1},\dots,x_{n}]$ will then be called a {\sl
cyclic nil-polynomial,} see also Table 1. In case $\FF=\RR$ the
quadratic form $f_{[2]}$ has type
$(\lceil{n\over2}\rceil,\lfloor{n\over2}\rfloor)$.

{{\def\quad{\hskip6pt}\klein
\setbox\strutbox=\hbox{\vrule height 10pt depth 6pt width 0pt}
\midinsert
$$\vbox{{\offinterlineskip\tabskip=0pt
\halign{\strut\vrule#&
 \quad\hfil$#$\hfil\quad& \vrule#&
 \quad\hfil$#$\hfil\quad& \vrule#&
 \quad\hfil$#$\hfil\quad& \vrule#&
 \quad\hfil$#$\hfil\quad& \vrule#&
 \quad \hfil$#$\hfil\quad& \vrule#\cr  
\noalign{\hrule}
&f_{[2]}&&f_{[3]}&&f_{[4]}&&f_{[5]}&&f_{[6]}&\cr
\noalign{\hrule}
&0&&0&&0&&0&&$0$&\cr
&x_{1}^{(2)}&&0&&0&&0&&$0$&\cr
&x_{1}x_{2}&&x_{1}^{(3)}&&0&&0&&0&\cr
&x_{1}x_{3}+x_{2}^{(2)}&&x_{1}^{(2)}\!x_{2}&&x_{1}^{(4)}&&0&&0&\cr
&x_{1}x_{4}+x_{2}x_{3}&&x_{1}x_{2}^{(2)}+x_{1}^{(2)}\!x_{3}&&
     x_{1}^{(3)}x_{2}&&x_{1}^{(5)}&&0&\cr
&x_{1}x_{5}+x_{2}x_{4}+x_{3}^{(2)}&&x_{1}x_{2}x_{3}+x_{1}^{(2)}
     \!x_{4}+x_{2}^{(3)}&&x_{1}^{(2)}\!x_{2}^{(2)}+x_{1}^{(3)}\!x_{3}&&
       x_{1}^{(4)}\!x_{2}&&x_{1}^{(6)}&\cr
\noalign{\hrule}
}}}
$$\vskip-8pt \centerline{\klein Table 1:~~Cyclic
nil-polynomials $f\in\FF[x_{1},\dots,x_{n-1}]$ for $1\le n\le6$, where
$y^{(k)}\!:=\qu{y^{k}\!}{(k!)}$} 
\endinsert}}\Formend

\medskip For graded PANAs $\5N=\bigoplus_{k>0}\5N_{k}$ with
annihilator $\5A=\5N_{d}$ we have for every $s\in\FF^{*}$ the algebra
automorphism
$$\theta_{s}:=\bigoplus_{k>0}s^{k}\id_{|\5N_{k}}\in\Aut(\5N)\,.\Leqno{KB}$$
As a consequence, if $t\in\FF^{*}$ admits a $d\,$\th root in $\FF$,
the pointings $\omega$ and $t\omega$ differ by an automorphism of
$\5N$. 

We mention that the PANA $\5N$ associated with the nil-polynomial
$f=\omega+\6q+\6c$ considered in Proposition \ruf{HU} also has a grading:
Indeed, put $\5N_{1}:=W_{1}$, $\5N_{2}:=W_{2}$ and $\5A:=\5N_{3}:=\FF$ with
products given by $xy:=x\cd y$ if $x,y\in W_{1}$ and
$xy:=\omega_{2}(x,y)$ if $x\in W_{1}$, $y\in W_{2}$. Using this we can
improve the second statement in \ruf{HU}.

\medskip PANAs $\5N$ admitting a grading enjoy a special property:
It is easy to see as a consequence of \Ruf{QL} that for every
associated nil-polynomial $f\in\FF[x_{1},\dots,x_{n}]$ there exist
linear forms $\lambda_{1},\dots,\lambda_{n}$ on $\FF^{n}$ with
$$f=\sum_{k=1}^{n}\lambda_{k}\df{x_{k}}\,.\Leqno{UR}$$

\Proposition{HW} With the notation of Proposition \ruf{HU} assume that
the cubic form $\6c$ on $W_{1}$ has the following property:\nline
$(*)$\hfill {\rm $z=0$ is the only element $z\in W_{1}$ with
$\6c(x+z)=\6c(x)$ for all $x\in W_{1}$.}\nline The graded PANA
$\5N=W_{1}\oplus W_{2}\oplus\FF$ with product $(x,y)\mapsto xy$
corresponding to the nil-polynomial $f=\6q+\6c$ on $W$ then
satisfies $\5N^{2}=W_{2}\oplus\FF$ as a consequence of
$(*)$. Furthermore, if $\;\tilde\6c$ is a second cubic form on $W_{1}$
with nil-polynomial $\tilde f=\6q+\tilde\6c$ and PANA
$\tilde\5N=W_{1}\oplus W_{2}\oplus\FF$ with appropriate product, the
following conditions are equivalent. \0 The nil-polynomials $f$,
$\tilde f$ are equivalent. \1 The algebras $\5N$, $\tilde\5N$ are
isomorphic as abstract algebras. \1 $\tilde f=f\circ g$ for some
$g\in\GL(W_{1})$.

\Proof Let $\omega:W\oplus\FF\to\FF$ be the canonical projection. Then
$\omega$ is a pointing for $\5N$. As in \Ruf{HO} define the
$\omega_{k}$ for the nil-polynomial $f=f_{[2]}+f_{[3]}$ on $W$. We have to
show that $\{x\cd y:x,y\in W_{1}\}$ spans $W_{2}$. If not, there
exists a vector $z\ne0$ in $W_{1}$ with $\omega_{2}(x\cd y,z)=0$ for
all $x,y\in W_{1}$. But then $\omega_{3}(x,y,z)=0$ for all $x,y\in
W_{1}$ implies $\6c(x+z)=\6c(x)$ for all $x\in W_{1}$, a
contradiction.

\noindent\TO12 This follows immediately from Definition
\ruf{HV}.

\noindent\To31 This follows immediately from the second claim
in Proposition \ruf{HU}.

\noindent\To23 Let $h:\5N\to\tilde\5N$ be an algebra isomorphism. Then
$h(W_{2}\oplus\FF)\subset\tilde\5N^{2}\subset (W_{2}\oplus\FF)$ as a
consequence of $(*)$, that is, there is a $g\in\GL(W_{1})$ with
$h(x)\equiv g(x)\mod\5N^{2}$ for all $x\in W_{1}$.  By the second
statement in Proposition \ruf{HU} we may assume $g=\id$ without loss
of generality. But then $\tilde\6c(x)=\6c(h(x))=\6c(x)$ for all $x\in
W_{1}$ implies $\tilde f=f$.\qed

\medskip Suppose that $W_{1}\cong\FF^{m}$ with coordinates
$(x_{1},\dots,x_{m})$ has dimension $m>0$ in Proposition \ruf{HW}. Then
$W\cong\FF^{2m}$ with coordinates
$(x_{1},\dots,x_{m},y_{1},\dots,y_{m})$ and we may assume
$\6q(x,y)=x_{1}y_{1}+\dots+x_{m}y_{m}$. As already mentioned, the
linear space $\6C$ of all cubic forms $\6c$ on $W_{1}$ has dimension
${m+2\choose3}$. The subset $\6C^{*}$ of all $\6c\in\6C$ satisfying
the condition $(*)$ in Proposition \ruf{HW} is Zariski open and dense
in $\6C$. The group $\GL(W_{1})$ acts on $\6C^{*}$ from the right and has
dimension $m^{2}$ over $\FF$. The difference of dimensions is
$m\choose3$.  But this number is also the cardinality of the subset
$J\subset\NN^{3}$, consisting of all triples $j=(j_{1},j_{2},j_{3})$
with $1\le j_{1}<j_{2}<j_{3}\le m$. Consider the affine map
$$\alpha:\FF^{J}\to\6C\,,\quad (t_{j})\mapsto \6c^{}_{0}+\sum_{j\in
J}t_{j}\6c_{j}\,,\Leqno{HX}$$ where
$\6c_{0}:=x_{1}^{3}+\dots+x_{m}^{3}\in\6C^{*}$
and $\6c_{j}:=x_{j_{1}}x_{j_{2}}x_{j_{3}}\in\6C$ for all $j\in J$.

In case $\FF=\RR$ or $\FF=\CC$, for a suitable neighbourhood $U$ of
$0\in\FF^{J}$ the map $\alpha:U\to\6C^{*}$ intersects all
$\GL(n,\FF)$-orbits in $\6C^{*}$ transversally. Indeed, since all
partial derivatives of $\6c_{0}$ are monomials containing a square,
the tangent space at $\6c_{0}$ of its $\GL(n,\FF)$-orbit is
transversal to the linear subspace $\langle\6c_{j}:j\in
J\rangle$ of $\6C$. In particular, in case $m\ge3$ there is a
family of dimension ${m\choose3}\ge1$ over $\FF$ $(=\RR\hbox{ or
}\CC)$ of pairwise different $\GL(n,\FF)$-orbits and thus of
non-equivalent nil-polynomials of degree $3$ on $W$. Notice that in
case $m=3$ the mapping $\alpha$ in \Ruf{HX} reduces to
$$\alpha:\FF\to\6C,\quad t\mapsto
x_{1}^{3}+x_{2}^{3}+x_{3}^{3}+tx_{1}x_{2}x_{3}\,.$$

We apply Proposition \ruf{KA} to the examples considered in
Proposition \ruf{HU}. Here $V=W_{1}\oplus W_{2}\oplus\FF$, $\6q$ is a
quadratic form on $W:=W_{1}\oplus W_{2}$ and $\6c$ is a cubic form on
$W_{1}$. With the extended nil-polynomial $f(x,y,t)=t+\6q(x+y)+\6c(x)$
on $W_{1}\oplus W_{2}\oplus\FF$ consider the hypersurface $$F:=\{z\in
V:f(z)=0\}$$ and identify $\FF$ with the line $\{0\}\oplus\FF$ in
$W\oplus\FF$. Then the affine group $\Aff(F)$ is transitive on $F$ by
Proposition \ruf{KA}, and every orbit in $V$ intersects the line
$\FF$. For every $s\in\FF^{*}$ the transformation $\theta_{s}$, see
\Ruf{KB}, satisfies $f\circ\theta_{s}=s^{3}h$ for every $s\in\FF^{*}$,
that is, $\theta_{s}\in\Aff(F)$. As a consequence, the number of
$\Aff(F)$-orbits in $V$ is bounded by the number of
$(\FF^{*})^{3}$-orbits in $\FF$. In particular, if
$(\FF^{*})^{3}=\FF^{*}$, then there are precisely two $\Aff(F)$-orbits
in $V$, namely $F$ and its complement. This situation occurs, for
instance, for $\FF=\RR$ and also for $\FF=\CC$. In any case, $F$ is
the only Zariski closed $\Aff(F)$-orbit in $V$, and every other orbit
is Zariski dense.

\bigskip \medskip \centerline{\bf Nil-polynomials of degree 4}

The method in Proposition \ruf{HU} can be generalized to get
nil-polynomials of higher degree, say of degree 4 for simplicity.
Throughout the subsection we use the notation \Ruf{KP}.

Let $W=W_{1}\oplus W_{2}\oplus W_{3}$ be a vector space with
$W_{1}=\FF^{n}$, $W_{2}=\FF^{m}$ and let $\6q$ be a fixed
non-degenerate quadratic form on $W$ in the following. Assume that
$W_{1}$, $W_{3}$ are totally isotropic and that $W_{1}\oplus W_{3}$,
$W_{2}$ are orthogonal with respect to $\6q$. Then $W$ has dimension
$2n+m$, and without loss of generality we assume that 
$$\6q(y)=\sum_{k=1}^{m}\epsilon_{k}y_{k}^{(2)}\Steil{for
suitable}\epsilon_{k}\in\FF^{*}\steil{and all}y\in W_{2}\,.$$ As
before let $\6C$ be the space of all cubic forms on $W$. Our aim is to
find cubic forms $c\in\6C_{\6q}$ that are the cubic part of a
nil-polynomial of degree 4.

\def\C{\6C'}

Denote by $\C$ the space of all cubic forms $\6c$ on $W_{1}\oplus
W_{2}$ such that $\6c(x+y)$ is quadratic in $x\in W_{1}$ and linear in
$y\in W_{2}$, or equivalently, which are of the form
$$\6c(x+y)={1\over2}\sum_{k=1}^{m}\sum_{i,j=1}^{n}\6c_{ijk}x_{i}x_{j}y_{k}
\Steil{for all}x\in W_{1},y\in W_{2}$$ with suitable coefficients
$\6c_{ijk}=\6c_{jik}\in\FF$.  Extending every $\6c\in\C$ trivially to
a cubic form on $W$ we consider $\C$ as subset of $\6C$.

For fixed $\6c\in\C$ the symmetric 2- and 3-linear forms
$\omega_{2},\omega_{3}$ on $W$ are defined by
$\omega_{2}(x,x)=2\6q(x)$ and $\omega_{3}(x,x,x)=6\6c(x)$ for all
$x\in W$. With the commutative product $x\cd y$ on $W$, see \Ruf{HY},
define in addition also the $k$-linear forms $\omega_{k}$ by \Ruf{HQ}
for all $k\ge4$. Then,
for every $x,y\in W_{1}$ the identity
$\omega_{2}(x\cd y,t)=\omega_{3}(x,y,t)=0$ for all $t\in W_{1}\oplus
W_{3}$ implies $x\cd y\in W_{2}$, that is $W_{1}\cd W_{1}\subset
W_{2}$. In the same way $\omega_{2}(x\cd y,t)=0$ for all $x\in W_{1}$,
$y\in W_{2}$ and $t\in W_{2}\oplus W_{3}$ implies $W_{1}\cd
W_{2}\subset W_{3}$. Also $W_{j}\cd W_{k}=0$ follows for all $j,k$
with $j+k\ge4$. Therefore $\6c$ belongs to $\6C_{\6q}$ if and only if
$(a\cd b)\cd c=a\cd(b\cd c)$ for all $a,b,c\in W_{1}$.

In terms of the standard basis $e_{1},\dots,e_{m}$ of $W_{2}=\FF^{m}$
we have
$$a\cd b=\sum_{k=1}^{m}\Big(\sum_{i,j=1}^{n}\epsilon_{k}^{-1}
\6c_{ijk}a_{i}b_{j}\Big) e_{k}\Steil{for all}a,b\in W_{1}$$ and thus
with $\Theta_{i,j,r,s}:=\sum_{k=1}^{m}\epsilon_{k}^{-1}
\6c_{ijk}\6c_{rsk}$ we get the identity
$$\eqalign{\omega_{2}\big((a\cd b)\cd c,t\big)&=\omega_{3}(a\cd
b,c,t)=\omega_{2}(a\cd b,c\cd
t)\cr&=\sum_{i,j,r,s=1}^{n}\!\!\Theta_{i,j,r,s}\,a_{i}
b_{j}c_{r}t_{s}\Steil{for all}a,b,c,t \in W_{1}\,.\cr}$$ Therefore, 
$$\6A:=\C\cap\6C_{\6q}=\{\6c\in\C:\Theta_{i,j,r,s}\steil{is symmetric
in}i,r\}\,.\Leqno{UF}$$ Notice that the condition in \Ruf{UF} implies
that $\Theta_{i,j,r,s}$ is symmetric in all indices. $\6A$ is a
rational subvariety of the linear space $\C$, it consists of all those
$\6c$ for which the corresponding product $x\cd y$ on $W$ is
associative.

\medskip For every $\6c\in\6A$ the corresponding
nil-polynomial $f$ on $W$ has the form
$$f=f_{[2]}+f_{[3]}+f_{[4]}\steil{with}f_{[2]}=\6q,\;f_{[3]}=\6c\hbox{~~and}$$
$$f_{[4]}(z)={1\over12}\6q(x\cd x)\steil{for all}z=(x,y,t)\in
W_{1}\oplus W_{2}\oplus W_{3}\,.$$ The group
$\Gamma:=\GL(W_{1})\times\O(\6q_{|W_{2}})\;\subset\;\GL(W_{1}\oplus
W_{2})$ acts on $\C$ by $\6c\mapsto\6c\circ\gamma^{-1}$ for every
$\gamma\in\Gamma$. Furthermore, $(g,h)\mapsto(g,h,(g^{\sharp})^{-1})$
embeds $\Gamma$ into $\O(\6q)$, compare the proof of Proposition
\ruf{HU}. As a consequence, the subvariety $\6A\subset\C$ is invariant
under $\Gamma$.

For every $\6c\in\6A$ the corresponding nil-polynomial comes from a
graded PANA with nil-index 4, provided $\6c\ne0$. Indeed, put
$W_{4}:=\FF$ and endow $W\oplus W_{4}$ with the product \Ruf{HR}.  It
is obvious that the linear span of $W_{1}\cd W_{1}$ in $W_{2}$ has
dimension $\le{n+1\choose2}$.

Let us consider the special case $n=2$ with $m={n+1\choose2}=3$ in
more detail. For simplicity we assume that for suitable coordinates
$(x_{1},x_{2})$ of $W_{1}$, $(y_{1},y_{2},y_{3})$ of $W_{2}$ and
$(z_{1},z_{2})$ of $W_{3}$ the quadratic form $\6q$ is given by
$$\6q=x_{1}z_{1}+x_{2}z_{2}+y_{1}^{(2)}+ y_{2}^{(2)}+\epsilon
y_{3}^{(2)}\Steil{for fixed}\epsilon\in\FF^{*}\Leqno{KU}$$ (in case
$\FF=\RR,\CC$ this is not a real restriction). For every $t\in\FF$
consider the cubic form
$$\6c_{t}:=(x_{1}^{(2)}+x_{2}^{(2)})y_{1}+x_{1}x_{2}y_{2}+
tx_{2}^{(2)}y_{3}$$ on $W_{1}\oplus W_{2}$.  A simple computation
reveals that every $\6c_{t}$ is contained in $\6A=\C\cap\6C_{q}$. The
corresponding nil-polynomial (depending on the choice of $\epsilon$)
then is
$$f_{t}=\6q+\6c_{t}+\6d_{t}\Steil{with}\6d_{t}:=x_{1}^{(4)}
+x_{1}^{(2)}x_{2}^{(2)}+(1+\epsilon^{-1}t^{2}) x_{2}^{(4)}\,.$$ In
addition we put
$$f_{\infty}:=\6q+x_{2}^{(2)}y_{3}+\epsilon^{-1}x_{2}^{(4)}$$ (a
smash product with the cyclic nil-polynomial of degree 4, see Table
1), where $\infty$ in the projective line
$\PP_{1}(\FF)=\FF\cup\{\infty\}$ is the {\sl point at
infinity}. Notice that for $t\in\FF^{*}$ the nil-polynomials $f_{t}$
and $\tilde f_{1/t}:=\6q+t^{-1}\6c_{t}+t^{-2}\6d_{t}$ are equivalent.

It is obvious that $f_{t}$ is equivalent to $f_{-t}$ for every
$t\in\PP_{1}(\FF)$.  Also, for every cubic term $\6c_{t}$ with
$t\in\FF^{*}$ the set $W_{1}\cd W_{1}$ spans $W_{2}$. For $t=0,\infty$
the linear span of $W_{1}\cd W_{1}$ in $W_{2}$ has dimension $2,1$
respectively.  For every $t\in\FF^{*}$ an invariant of $\6d_{t}$ is
the number $\pphi(t):=\qu{g_{2}(\6d_{t})^{3}}{g_{3}(\6d_{t})^{2}}=
\epsilon^{2}t^{-4}(4+\epsilon^{-1}t^{2})^{3}\in\FF$, where
$g_{2},g_{3}$ are the classical invariants of binary quartics, compare
\Lit{MUKA} \p27.  Since every fiber of $\pphi:\FF^{*}\to\FF$ contains
at most 6 elements we conclude

\Proposition{} For every Field $\FF$ and every fixed
$\epsilon\in\FF^{*}$ the set of all equivalence classes given by all
nil-polynomials $f_{t}$, $t\in\FF$, has the same cardinality as
$\FF$ and, in particular, is infinite.  \Formend

{\medskip\noindent\bf Remarks} 1. In case $\FF=\QQ$ is the rational
field there are infinitely many choices of $\epsilon\in\QQ^{*}$
leading to pairwise non-equivalent quadratic forms $\6q$ in \Ruf{KU}.
For each such choice there is an infinite number of pairwise
non-equivalent nil-polynomials $f_{t}$ of degree 4 over $\QQ$.

\noindent 2. In case $\FF=\RR$ is the real field there are essentially
the two choices $\epsilon=\pm1$. In case $\epsilon=1$ the form $\6q$
has type $(5,2)$ and all nil-polynomials $f_{t}$ with $0\le
t\le\sqrt8$ are pairwise non-equivalent.  In case $\epsilon=-1$ the
form $\6q$ has type $(4,3)$ and all $f_{t}$ with $0\le t\le\infty$ are
pairwise non-equivalent.

\noindent 3. Nil-polynomials of degree $\ge5$ can be
constructed just as in the case of degrees 3 and 4 as before. As an
example we briefly touch the case of degree 5: Fix a vector
space $W$ of finite dimension over $\FF$ together with a
non-degenerate quadratic form $\6q$ on $W$. Assume furthermore that
there is a direct sum decomposition $W=W_{1}\oplus W_{2}\oplus
W_{3}\oplus W_{4}$ into totally isotropic subspaces such that
$W_{1}\oplus W_{4}$ and $W_{2}\oplus W_{3}$ are orthogonal. Then
consider a cubic form $\6c$ on $W$ that can be written as a sum
$\6c=\6c'+\6c''$ of cubic forms with the following properties: $\6c'$
is a cubic form on $W_{1}\oplus W_{3}$ (trivially extended to $W$)
that is linear in the variables of $W_{3}$ while $\6c''$ is a cubic
form on $W_{1}\oplus W_{2}$ that is linear in the variables of
$W_{1}$. Denote by $x\cd y$ the commutative product on $W$ determined
by $\6q$ and $\6c$. If we put $W_{k}:=0$ for all $k>4$ we have
$W_{j}\cd W_{k}\subset W_{j+k}$ for all $j,k$. Therefore,
$\6c\in\6C_{\6q}$ if and only if the product $x\cd y$ on $W$ is
associative, see \Ruf{QH} for the notation. This is true without any
assumption if $\6c=\6c'$ or $\6c=\6c''$. But $W_{1}\cd W_{1}=0$ in the
first and $W_{2}\cd W_{2}=0$ in the latter case. On the other hand,
the nil-polynomial associated to $\6c\in\6C_{\6q}$ has degree 5 if
$W_{2}\cd W_{2}$ spans $W_{4}$ and $W_{1}\cd W_{1}\ne0$.

\bigskip \medskip \centerline{\bf Affine homogeneity}

\medskip In this subsection let $\FF$ be either $\RR$ or $\CC$.  Also
let $\5N\ne0$ be an arbitrary commutative associative nilpotent
algebra of finite dimension over $\FF$. In addition we assume that
there exists a $\ZZ$-gradation
$$\5N=\bigoplus_{k>0}\5N_{k}\,,\quad\5N_{j}\5N_{k}\subset\5N_{j+k}\,.$$
Let $d:=\max\{k:\5N_{k}\ne0\}$ and denote by $\pi:\5N\to\5N_{d}$ the
canonical projection with kernel $\5K:=\bigoplus_{k<d}\5N_{k}$. We do
not require that $\5N_{k}\ne0$ for all $1\le k\le d$ nor that
$\5N_{d}$ is the annihilator or has dimension 1. Extending $\pi$
linearly to $\5N^{0}$ by $\pi(\One)=0$ we have the polynomial map
$\6f:=\pi\circ\exp:\5N\to\5N_{d}$. The submanifold $F:=\6f^{-1}(0)$
then is the graph of a polynomial map $\5K\to\5N_{d}$ and
$\5K=\T_{0}F$ is the tangent space to $F$ at the origin. We are
interested in the affine group $\Aff(F)=\{g\in\Aff(\5N):g(F)=F\}$ and
its subgroup $\6A=\6A(\6f):=\{g\in\Aff(\5N):\6f\circ g=\6f\}$.

Every point $x\in\5N$ has a unique representation $x=x_{1}+\dots+x_{d}$
with $x_{k}\in\5N_{k}$. Consider on $\5N$ the linear span $\7a$ of all
nilpotent affine vector fields of the form
$$(d-j)\alpha_{j}\!\dd{x_{j}}-\;\sum_{k=1}^{d-j}k\alpha_{j}x_{k}\dd{x_{j+k}}\Steil{with}1\le j<
d\steil{and}\alpha_{j}\in\5N_{j}\,.$$ 

{\klein\openup-2pt As an example, in case
$\5N\cong\FF^{4}$ is the cyclic PANA of nil-index 4, see Table 1, we
have $d=4$,\nline
$\,\6f(x)=x_{4}+x_{1}x_{3}+x_{2}^{(2)}+x_{1}^{(2)}\!x_{2}+x_{1}^{(4)}$,
and $\7a$ is the linear span of the vector fields
$$\eqalign{
3\dd{x_{1}}-x_{1}\dd{x_{2}}-2x_{2}\dd{x_{3}}&-3x_{3}\dd{x_{4}}\cr
2\dd{x_{2}}-\phantom{2}x_{1}\dd{x_{3}}&-2x_{2}\dd{x_{4}}\cr
\dd{x_{3}}&-\phantom{2}x_{1}\dd{x_{4}}\;.\cr }$$\par} With some
computation we get: 

\Lemma{GR} $\7a$ is a nilpotent Lie algebra and
the evaluation map $\epsilon_{a}:\7a\to\5N$, $\;\xi\mapsto\xi_{a}$, is
injective for every $a\in\5N$. In particular, all orbits in $\5N$ of
the nilpotent subgroup $\exp(\7a)\subset\Aff(\5N)$ have the same
dimension.

\Proposition{} $\7a\!\6f=0$, that is, $\exp(\7a)\subset\6A(f)$. In
particular, $\6A(\6f)$ acts transitively on every $c+F=\6f^{-1}(c)$,
$\,c\in\5N_{d}$.
 
\Proof Put $\xi:=(d-j)\alpha\!\dd{x_{j}}-\;\sum_{k=1}^{d-j}
k\alpha x_{k}\dd{x_{j+k}}$ for fixed $1\le j< d$ and
$\alpha\in\5N_{j}$. Then
$$\xi\6f=\sum c_{\nu}\alpha\,x_{1}^{(\nu_{1})}x_{2}^{(\nu_{2})}\cdots
x_{d}^{(\nu_{d})}\,,$$ where the sum is taken over all multi indices
$\nu\in\NN^{d}$ with $\nu_{1}+2\nu_{2}+\dots+d\nu_{d}=d-j$ and
$c_{\nu}$ certain rational coefficients. Now fix such a multi index
$\nu$. For simpler notation we put $x^{(-1)}:=0$ for every $x\in\5N$.
Then we have {\klein
$$\eqalign{c_{\nu}^{}\alpha\,x_{1}^{(\nu_{1})}x_{2}^{(\nu_{2})}\cdots
x_{d}^{(\nu_{d})}&=(d-j)\alpha\!\dd{x_{j}}\Big(x_{1}^{(\nu_{1})}\cdots
x_{j}^{(\nu_{j}+1)}\cdots x_{d}^{(\nu_{d})}\Big)\,-\cr
&\phantom{XXXX}\sum_{k=1}^{d-j}k\alpha
x_{k}\dd{x_{j+k}}\Big(x_{1}^{(\nu_{1})}\cdots
x_{k}^{(\nu_{k}-1)}\cdots x_{j+k}^{(\nu_{j+k}+1)}\cdots
x_{d}^{(\nu_{d})}\Big)\cr
&=\big(d-j-\sum_{k=1}^{d-j}k\nu_{k}\big)\alpha\,x_{1}^{(\nu_{1})}
x_{2}^{(\nu_{2})}\cdots x_{d}^{(\nu_{d})}\;=\;0 }$$} since $\nu_{k}^{}=0$
for $k>d-j$.\qed

\bigskip Next we specialize to the case where $\5N_{d}$ has dimension
1, that is, $F$ is a hypersurface in $\5N$ (we still do not require
that $\5N_{d}$ is the annihilator of $\5N$ although contained in it).
For every $s\in\FF^{*}$ we have the semi-simple linear transformation
$$\theta_{s}:=\bigoplus_{k>0}s^{k}\id_{|\5N_{k}}\in\GL(F)$$ satisfying
$\6f\circ\theta_{s}=s^{d}\6f$. As a consequence we have that the group
$\Aff(F)$ has at most 3 orbits in $\5N$. In case $d$ odd or $\FF=\CC$
this group has only two orbits in $\5N$, the closed hypersurface $F$
and the open complement $\5N\backslash F$. As in the subsection
\ruf{complex} we get in case $\FF=\CC$ further affinely homogeneous
real surfaces.

\bigskip

\vskip7mm {\gross\noindent References} \medskip {\klein\openup-.5pt
\parindent 15pt\advance\parskip-1pt

\def\Springer{Ber\-lin-Hei\-del\-berg-New York: Sprin\-ger~}

\Ref{BAJT}Baouendi, M.S., Jacobowitz, H., Treves F.: On the analyticity of CR mappings. Ann. of Math. {\bf 122} (1985), 365-400.
\Ref{BERO}Baouendi, M.S., Ebenfelt, P., Rothschild, L.P.: {\sl Real Submanifolds in Complex Spaces and Their Mappings}. Princeton Math. Series {\bf 47}, Princeton Univ. Press, 1998.
\Ref{CHMO}Chern, S.S., Moser, J.K.: Real hypersurfaces in complex manifolds. Acta. Math. {\bf 133} (1974), 219-271.
\Ref{COUR}Courter, R.C.: Maximal commutative subalgebras of $K_{n}$  at exponent three. Linear Algebra and Appl. {\bf 6} (1973), 1-11.
\Ref{DAYA}Dadok, J., Yang, P.: Automorphisms of tube domains and spherical hypersurfaces. Amer. J. Math. {\bf 107} (1985), 999-1013.
\Ref{EAST}Eastwood, M.G.: Moduli of isolated hypersurface singularities. Asian J. Math. {\bf 8} (2004), 305-314.
\Ref{EAEZ}Eastwood, M.G., Ezhov, V.V.: On affine normal forms and a classification of homogeneous surfaces in affine three-space. Geom. Dedicata  {\bf 77} (1999), 11-69.
\Ref{FELS}Fels, G.: Locally homogeneous finitely nondegenerate CR-manifolds. Math. Res. Lett. {\bf 14} (2007), 693-922.
\Ref{FKAU}Fels, G., Kaup, W.: Classification of Levi degenerate homogeneous CR-manifolds of dimension 5. Acta Math. {\bf 201} (2008), 1-82.
\Ref{FKAP}Fels, G., Kaup, W.: Local tube realizations of CR-manifolds
and maximal abelian subalgebras. Ann. Scuola Norm. Sup. Pisa. To appear
\Ref{ISMI}Isaev, A.V., Mishchenko, M.A.: Classification of spherical tube hypersurfaces that have one minus in the Levi signature form. Math. USSR-Izv. {\bf 33} (1989), 441-472.
\Ref{ISAE}Isaev, A.V.: Classification of spherical tube hypersurfaces that have two minuses in the Levi signature form. Math. Notes {\bf 46} (1989), 517-523.
\Ref{ISAV}Isaev, A.V.: Global properties of spherical tube hypersurfaces. Indiana Univ. Math. {\bf 42} (1993), 179-213.
\Ref{ISEV}Isaev, A.V.: {\sl Rigid and Tube Spherical Hypersurfaces.} Book manuscript in preparation. 
\Ref{KAPP}Kaup, W.: On the holomorphic structure of $G$-orbits in compact hermitian symmetric spaces. Math. Z. {\bf 249} (2005), 797-816.
\Ref{KAZT}Kaup, W., Zaitsev, D.: On local CR-transformations of Levi degenerate group orbits in compact Hermitian symmetric spaces. J. Eur. Math. Soc. {\bf 8} (2006), 465-490.
\Ref{KNAP}Knapp, Anthony W.: Lie groups beyond an introduction. 2nd ed. Progress in Mathematics {\bf 140},  Birkh\"auser, Boston 2002.
\Ref{LAFF}Laffay, Th.J.: The Minimal Dimension of Maximal Commutative Subalgebras of Full Matrix Algebras. Linear Algebra and Appl. {\bf 71} (1985), 199-212.
\Ref{MUKA}Mukai, S.: An introduction to invariants and moduli. Cambridge Univ. Press 2003
\Ref{TANA}Tanaka, N.: On the pseudo-conformal gem,ometry of hypersurfaces of the space of $n$ complex variables. J. Math. Soc. Japan {\bf 14} (1962), 397-429.
\Ref{WARN}Warner, G.: Harmonic analysis on semi-simple Lie groups I. Die Grundlehren der mathematischen Wissenschaften. Band {\bf188}. \Springer  1972.

\bigskip
}

\smallskip\openup-4pt\parindent0pt
\hbox to 3cm{\hrulefill}
G. Fels\par e-mail: {\tt gfels@uni-tuebingen.de\par\vskip8pt}
W. Kaup\par e-mail: {\tt kaup@uni-tuebingen.de\par\vskip8pt}

Mathematisches Institut, Universit\"at
T\"ubingen,\par Auf der Morgenstelle 10,\par  72076 T\"ubingen,
Germany\par

\closeout\aux\bye